\documentclass[
	english
]{scrartcl}

\usepackage[T1]{fontenc}
\usepackage[utf8]{inputenc}
\usepackage[english]{babel}
\usepackage[inline,shortlabels]{enumitem}
\usepackage{amsmath}
\usepackage{amssymb}
\usepackage{amsthm}
\usepackage[noadjust]{cite}
\usepackage[a4paper,total={6in,9in}]{geometry}
\usepackage[colorlinks,citecolor=blue]{hyperref}
\usepackage[capitalize,noabbrev,nameinlink]{cleveref}
\usepackage{xcolor}
\usepackage{longtable}
\usepackage{tikz,pgfplots,pgfplotstable}
\usetikzlibrary{patterns}
\usepackage{algorithmic}
\usepackage{marginnote}
\usepackage[labelformat=simple]{subcaption}

\setlist{font=\normalfont,topsep=1ex,parsep=0ex}

\setlist[enumerate]{label=(\alph*)}

\usepackage{changes}

\numberwithin{equation}{section}
\numberwithin{table}{section}
\numberwithin{figure}{section}

\usepackage[figure]{hypcap}
\usepackage{graphicx}
\graphicspath{{./img/}}

\usepackage{preamble}

\crefname{figure}{Figure}{Figures}
\crefname{table}{Table}{Tables}
\crefname{assumption}{Assumption}{Assumptions}
\Crefname{ALC@unique}{Step}{Steps}

\newlist{alglist}{enumerate}{1}
\setlist[alglist]{topsep=1ex,parsep=0ex,leftmargin=*,label=\textbf{Step~\arabic*.}}

\let\eqref\labelcref

\allowdisplaybreaks

\newcommand\norm[1]{\left\Vert#1\right\Vert}
\newcommand\nnorm[1]{\Vert#1\Vert}

\newcommand\N{\mathbb{N}}
\newcommand\R{\mathbb{R}}

\newcommand\tto{\rightrightarrows}

\newcommand{\cl}{\operatorname{cl}}

\newcommand{\intr}{\operatorname{int}}

\newcommand{\conv}{\operatorname{conv}}
\newcommand{\cone}{\operatorname{cone}}

\newcommand{\dom}{\operatorname{dom}}

\newcommand{\gph}{\operatorname{gph}}

\DeclareMathOperator*{\argmin}{\operatorname{argmin}}

\newcommand{\eff}[3]{\mathcal{E}( #1 , #2 , #3 )}			        
\newcommand{\nd}[3]{\mathcal{N}( #1 , #2 , #3 )}		            	
\newcommand{\weff}[3]{\mathcal{E}_{\textup{w}}( #1 , #2 , #3 )}	    
\newcommand{\wnd}[3]{\mathcal{N}_{\textup{w}}( #1 , #2 , #3 )}		

\DeclareMathAlphabet{\mathpzc}{OT1}{pzc}{m}{it}

\newtheorem{theorem}{Theorem}[section]
\newtheorem{lemma}[theorem]{Lemma}
\newtheorem{proposition}[theorem]{Proposition}
\newtheorem{assumption}[theorem]{Assumption}

\newtheorem{corollary}[theorem]{Corollary}
\newtheorem{remark}[theorem]{Remark}
\newtheorem{definition}[theorem]{Definition}
\newtheorem{example}[theorem]{Example}

\definecolor{mygreen}{rgb}{0.0,0.7,0.0}
\definecolor{mybrown}{rgb}{0.5,0.5,0.0}

\begin{document}

\title{Notes on the value function approach to multiobjective bilevel optimization}
\author{%
	Daniel Hoff%
	\footnote{%
		Ilmenau University of Technology,
		Institute of Mathematics, 
		98684 Ilmenau, 
		Germany, \email{daniel.hoff@tu-ilmenau.de},
		\orcid{0009-0001-2725-2692}
		}%
	\and
	Patrick Mehlitz%
	\footnote{%
		University of Duisburg-Essen,
		Faculty of Mathematics,
		45127 Essen,
		Germany,
		\email{patrick.mehlitz@uni-due.de},
		\orcid{0000-0002-9355-850X}%
		}%
	}

\publishers{}
\maketitle

\begin{abstract}
	This paper is concerned with the value function approach to multiobjective bilevel optimization
	which exploits a lower-level frontier-type mapping in order to replace the hierarchical model of
	two interdependent multiobjective optimization problems	by a single-level multiobjective
	optimization problem.
	As a starting point, different value-function-type reformulations are suggested and their
	relations are discussed.
	Here, we focus on the situations where the lower-level problem is solved up to efficiency
	or weak efficiency, and an intermediate solution concept is suggested as well.
	We study the graph-closedness of the associated efficiency-type and frontier-type mappings.
	These findings are then used for two purposes.
	First, we investigate existence results in multiobjective bilevel optimization.
	Second, for the derivation of necessary optimality conditions via the value function approach,
	it is inherent to differentiate frontier-type mappings in a generalized way.
	Here, we are concerned with the computation of upper coderivative estimates for the frontier-type
	mapping associated with the setting where the lower-level problem is solved up to
	weak efficiency. We proceed in two ways, relying, on the one hand, on a weak domination property
	and, on the other hand, on a scalarization approach.
	Illustrative examples visualize our findings and some flaws in the related literature.
\end{abstract}


\begin{msc}	
	\mscLink{49J52}, \mscLink{49J53}, \mscLink{90C29}, \mscLink{90C30}
\end{msc}

\section{Introduction}\label{sec:introduction}

In this note, we are concerned with multiobjective bilevel optimization problems of type
\begin{equation}\label{eq:BPP}\tag{BPP}
	\min\limits_{x,y}{}_{\mathcal{K}}\{F(x,y)\,|\,x\in X,\,y\in\widehat{\Psi}(x)\}
\end{equation}
where $F\colon\R^n\times\R^m\to\R^p$ is a continuous vector-valued function,
$\mathcal{K}\subset\R^p$ is a nonempty, closed, convex, pointed cone with nonempty interior,
$X\subset\R^n$ is a nonempty, closed set, and
$\widehat{\Psi}\colon\R^n\tto\R^m$ is a suitable solution mapping 
associated with the parametric multiobjective optimization problem 
\begin{equation}\label{eq:parametric_problem}\tag{P$(x)$}
	\min\limits_y{}_{\mathcal C}\{f(x,y)\,|\,y\in\Gamma(x)\}.
\end{equation}
Above, $f\colon\R^n\times\R^m\to\R^q$ is a continuous vector-valued function, 
$\mathcal C\subset\R^q$ is a nonempty, closed, convex, pointed cone with nonempty interior, and
$\Gamma\colon\R^n\tto\R^m$
is a set-valued mapping with a closed graph.
In our model problem, $\mathcal K$ and $\mathcal C$ play the role of ordering cones, i.e.,
they specify what ``minimization''
means in \eqref{eq:BPP} and \eqref{eq:parametric_problem},
up to different notions of optimality/efficiency in multiobjective optimization. 
Indeed, even for a fixed ordering cone, there exist diverse concepts of efficiency which aim to
characterize the valuable points in the feasible set when faced with multiple objective functions.
In this paper, we focus our attention on so-called efficiency and weak efficiency.
Exemplary, we refer the interested reader to the monographs \cite{Ehrgott2005,Jahn2011} for a detailed
introduction to multiobjective optimization and to \cref{sec:vector_optimization} where
we summarize some essential foundations of this topic. 

In the case where $p:=q:=1$ and $\mathcal K:=\mathcal C:=\R_+$, \eqref{eq:BPP} reduces to
a so-called standard bilevel optimization problem which is used to model situations of
interdependent decision making between two parties that cooperate with each other.
Let us note that this model is closely related to the so-called optimistic approach to bilevel
optimization as shown in \cite[Proposition~6.9]{DempeMordukhovichZemkoho2012}.
In the literature, \eqref{eq:BPP} and \eqref{eq:parametric_problem} are referred to
as the upper- and lower-level problem, respectively,
	and we will stick to this terminology even in the more general situation considered here.
Due to numerous underlying applications from, e.g., 
finance, economics, data as well as natural sciences, and energy management,
bilevel optimization can be found among the most active fields of research within the
optimization community, see the survey papers \cite{Dempe2003,Dempe2020}. 

For an introduction to bilevel optimization, we refer the
interested reader to the monographs 
\cite{Bard1998,Dempe2002,DempeKalashnikovPerezValdesKalashnykova2015,ShimizuIshizukaBard1997}.
We emphasize that bilevel optimization problems are, in general, highly irregular and nonconvex
optimization problems which are only implicitly given as a closed-form representation of
the solution mapping $\widehat\Psi$ associated with \eqref{eq:parametric_problem} 
is often not available. This causes the model \eqref{eq:BPP} to be rather challenging
from a theoretical as well as numerical perspective.

Since the set of efficient or weakly efficient points of a multiobjective optimization problem
is, usually, not a singleton, the minimization of a scalar function over the set of
efficient or weakly efficient points, in order to reduce the number of reasonable and 
economically utilizable points, is a 
	closely related
topic, see e.g.\ 
\cite{Bolintineanu1993,Bolintineanu1993b,BonnelKaya2010,HorstThoai1999},
and such problems already possess a hierarchical structure.
So-called semivectorial bilevel optimization problems, where only $p:=1$ and $\mathcal K:=\R_+$
are demanded in \eqref{eq:BPP}, i.e., only the underlying parametric optimization problem
\eqref{eq:parametric_problem} possesses multiple objective functions, 
provide a much more general model paradigm and have been investigated, e.g., in
\cite{Bonnel2006,BonnelMorgan2006,DempeGadhiZemkoho2013,DempeMehlitz2019,LiuWanChenWang2014,Zemkoho2016}.
The even more general situation where the objective functions of both decision makers in \eqref{eq:BPP} 
are allowed to be vector-valued has been considered, e.g., in
\cite{DandurandGuarneriFadelWiecek2014,Eichfelder2010,Eichfelder2020,GebhardtJahn2009,LafhimZemkoho2022,RuuskaMiettinenWiecek2012}.
For formal completeness, let us also mention that the setting where $q:=1$ and $\mathcal C:=\R_+$
hold while the objective function of \eqref{eq:BPP} is vector-valued is also reasonable,
see e.g.\ \cite{GadhiDempe2012,Ye2011,Yin2002}. 
It is well known from the literature that scalarization techniques can be used to transfer
multiobjective optimization problems into scalar ones, see e.g.\ \cite[Section~5]{Jahn2011}.
Hence, whenever the lower-level problem \eqref{eq:parametric_problem} possesses multiple
objective functions, it seems to be 
an obvious idea
to replace it with a scalarized counterpart and
to treat the scalarization parameters as additional upper-level variables.
This procedure, however, comes with a price. 
It has been shown in the setting of semivectorial bilevel optimization that this reformulation
induces additional local minimizers which might be attractive for solution algorithms,
see \cite{DempeMehlitz2019}.
Even worse, the more detailed study \cite{BenkoMehlitz2021} revealed that this transformation also
induces additional stationary points while the constraint qualifications needed to tackle the
reformulated problem might be stronger than the ones applicable to the original model.
These observations are likely to
	concern
multiobjective bilevel optimization problems, too.
Hence, one should not tackle \eqref{eq:parametric_problem} with scalarization approaches 
for the theoretical and numerical treatment of \eqref{eq:BPP} in a lightheaded way.

From the viewpoint of applications, the investigation of multiobjective bilevel optimization
problems is rather relevant. 
Let us point the reader to three interesting model problems arising in practice.
First, for example, in \cite{Yin2002}, transportation planning systems are considered 
in which a system manager wants to determine road tolls for a transportation network 
by deciding on the charges on each edge of the network. 
Thereby, the leader wants to minimize the total network cost and 
to maximize the total revenue as well as the consumers' surplus at the same time. 
For the lower-level problem, the travelers choose their routes in the network 
by the traffic volume. This circumstance is modeled as an equilibrium traffic assignment problem
in which each user chooses the most convenient path selfishly.	
Second, \cite{LuoLiuLiu2021} is concerned with standalone hybrid renewable energy systems 
which can reliably meet the energy demand of users in remote areas 
and reduce the intermittency of different renewable energy sources. 
At the upper-level stage, policy makers like a government give incentive subsidies 
to investors to fill in the gap between the market price and the price charged to consumers.
Thereby, the decision maker is interested to minimize the negative impact 
of hybrid renewable energy systems on the environment measured 
by the emitted annualized carbon dioxide and, conflicting to this, the subsidies 
borne by the government. 
The investor's objective for the lower-level problem is to reduce the annualized cost 
of the energy system project which comprise the annualized capital cost, 
annual fuel cost, annual system maintenance cost, and subsidies. 
In the constraints, technical components of several energy technologies are considered.  
Third, in \cite{MaZhangFengLevLi2021}, multiobjective optimization problems are used 
to model waste management with obnoxious effects. 
The government decides at the upper level about the location and scale of waste 
recycling centers to minimize the total economic cost and, simultaneously, 
the abominable effects caused by recycling centers under given capacity constraints. 
Based on the locations of the recycling centers, a sanitation company determines 
waste collection routing plans to minimize the vehicle traveling and fixed cost for which, exemplary, 
flow conservation for delivery, locations of the underlying routing graph, 
load of the vehicle on its edges, and the recycling center capacities 
are considered as constraints.

In this note, we are concerned with the so-called value function approach to bilevel optimization,
which dates back to \cite{Outrata1990} in the scalar case.
Let us, for a moment, assume that \eqref{eq:parametric_problem} is a scalar problem, i.e., $q:=1$ and $\mathcal C:=\R_+$.
With the aid of the so-called lower-level optimal value function $\varphi\colon\R^n\to\overline{\R}$ given by
\begin{equation}\label{eq:optimal_value_function}
	\forall x\in\R^n\colon\quad
	\varphi(x):=\inf\limits_y\{f(x,y)\,|\,y\in\Gamma(x)\},
\end{equation}
where $\overline{\R} := \R \cup \{ -\infty,\infty\}$ and  $\inf \emptyset  := \infty$, 
one can equivalently reformulate \eqref{eq:BPP} by means of
\begin{equation}\label{eq:VFR_classic}
	\min\limits_{x,y}{}_{\mathcal{K}}\{F(x,y)\,|\,x\in X,\,f(x,y)\leq\varphi(x),\,y\in\Gamma(x)\},
\end{equation}
as the feasible sets of \eqref{eq:BPP} and \eqref{eq:VFR_classic} coincide.
Noting that $\varphi$ is only implicitly available while being inherently nonsmooth,
and observing that classical constraint qualifications from nonsmooth programming
are likely to fail at all feasible points of \eqref{eq:VFR_classic}, see
\cite[Proposition~3.2]{YeZhu1995}, \eqref{eq:VFR_classic} is still a rather challenging
problem. Nevertheless, the so-called value function reformulation \eqref{eq:VFR_classic}
has turned out to be useful for the construction of necessary optimality conditions
and solution algorithms for \eqref{eq:BPP} in the setting where, additionally, $p:=1$ and $\mathcal K:=\R_+$,
see e.g.\ \cite{DempeDuttaMordukhovich2007,DempeFranke2016,DempeZemkoho2013,MordukhovichNamPhan2012,YeZhu1995} and the references therein.

Based on its successful application in the scalar situation, the paper \cite{LafhimZemkoho2022}
suggests an extension of the value function approach to multiobjective bilevel optimization.
In order to recover the approach of \cite{LafhimZemkoho2022}, 
let us introduce the set-valued mapping $\widehat{\Phi}\colon\R^n\tto\R^q$ by means of
\[
	\forall x\in\R^n\colon\quad
	\widehat{\Phi}(x):=f(x,\widehat{\Psi}(x)):=\{f(x,y)\,|\,y\in\widehat{\Psi}(x)\},
\]
and observe that \eqref{eq:BPP} is closely related to
\begin{equation}\label{eq:value_function_reformulation}\tag{VFR}
	\min\limits_{x,y}{}_{\mathcal{K}}\{F(x,y)\,|\,x\in X,\,f(x,y)\in\widehat{\Phi}(x),\,y\in\Gamma(x)\}
\end{equation}
which, again, we will refer to as value function reformulation of \eqref{eq:BPP}.
Let us note at this point that \eqref{eq:BPP} and \eqref{eq:value_function_reformulation} do not need to
be equivalent anymore in general, see \cref{ex:order_of_definitions_for_intermediate_approach_is_important}.
In \cite{LafhimZemkoho2022}, the authors exploit the model problem \eqref{eq:value_function_reformulation} for
the derivation of necessary optimality conditions for the associated model \eqref{eq:BPP} where the
lower-level parametric optimization problem is solved up to efficient points.
More precisely, the authors of \cite{LafhimZemkoho2022} transfer the classical arguments from \cite{YeZhu1995},
which 
apply to
 bilevel optimization problems with scalar upper- and lower-level objective functions,
to the multiobjective situation. Among others, this makes the (generalized) differentiation of the so-called
frontier mapping $\widehat\Phi$ from above necessary. Therefore, Mordukhovich's coderivative construction
is exploited, see e.g.\ \cite{Mordukhovich2006,RockafellarWets2009}. 
Let us mention that the computation of coderivatives associated with frontier maps has already been studied 
in the literature before, see e.g.\ \cite{HuyMordukhovichYao2008,LiXue2014,XueLiLiaoYao2011}.

Now, we would like to mention an important and crucial observation which seemingly has been 
	overlooked
in most of the aforementioned literature dealing with multiobjective bilevel optimization.
If the lower-level problem \eqref{eq:parametric_problem} is solved up to efficiency, then the associated
solution mapping is likely to possess a graph which is not closed - even under comparatively strict assumptions.
This non-closedness also 
	concerns
the associated frontier mapping which is used in the definition of the
associated value function reformulation. 
Let us point out two immediate disadvantages which result from this non-closedness.
First, it will be difficult to guarantee that \eqref{eq:BPP} possesses efficient or at least weakly efficient
points in this case.
Second, for the computation of coderivatives, the underlying set-valued mapping needs to possess,
at least locally around the reference point, a closed graph, so, in turn, 
this concept formally cannot be applied to the frontier mapping of interest in general.

In this note, we 
tackle
both of the aforementioned issues and, along the way, comment on related peculiarities
of multiobjective bilevel optimization. 
Throughout the paper, several examples are presented to illustrate our arguments.
Furthermore, we compare our results with related findings 
applying to
\eqref{eq:VFR_classic} in the scalar case.
To start, let us mention that whenever \eqref{eq:parametric_problem} is solved up to weak efficiency only,
then the associated solution map is likely to possess a closed graph under mild assumptions, and this
also extends to the associated frontier mapping.
This can be distilled e.g.\ from \cite[Section~3.5]{GoepfertRiahiTammerZalinescu2003}, but we provide
a self-contained proof here for the purpose of completeness. 
In \cref{sec:model_problem} of the paper, 
we are not only concerned with some formal properties of the value function reformulation,
but also with the derivation of existence results for \eqref{eq:BPP} which is a closely related issue.
As already mentioned, this is reasonable under mild assumptions if \eqref{eq:parametric_problem} is solved
up to weak efficiency. 
We also show that the efficiency and associated frontier mapping 
of a fully linear parametric multiobjective optimization problem possess closed graphs,
giving rise to yet another criterion for the existence of efficient and weakly efficient points for \eqref{eq:BPP}.
Furthermore, we suggest an intermediate solution approach between efficiency and weak efficiency 
applying to 
\eqref{eq:parametric_problem} which
naturally comes along with a closed graph of the associated mappings $\widehat\Psi$ as well as $\widehat\Phi$
and, thus, the superordinate problem \eqref{eq:BPP} is also likely to possess efficient and weakly efficient points.
\Cref{sec:coderivative_calculus} is dedicated to the computation of upper coderivative estimates for
the weak frontier mapping, i.e., the frontier mapping associated with \eqref{eq:parametric_problem} when solved
up to weak efficiency. First, we follow the approach promoted in \cite{HuyMordukhovichYao2008},
which relies on the so-called domination property, see e.g.\ \cite{Henig9867,Luc1984}, and point out
some associated advantages and disadvantages. 
Furthermore, we illustrate a crucial error in \cite[Theorem~4.1]{LafhimZemkoho2022}, that is concerned with the
computation of an upper estimate for the coderivative of the efficiency-related frontier mapping of \eqref{eq:parametric_problem}.
Unluckily, this shows that the results obtained in \cite[Section~4]{LafhimZemkoho2022} are not reliable.
A second approach to coderivative estimates we discuss here is based on a weighted-sum scalarization of the
problem \eqref{eq:parametric_problem} in the convex situation. This idea is motivated by considerations from \cite{Zemkoho2016}.

The remainder of the paper is organized as follows.
In \cref{sec:preliminaries}, we comment on the basic notation used in this paper, 
as well as some preliminaries from set-valued analysis and generalized differentiation. 
Furthermore, some foundations of multiobjective optimization are presented.
\Cref{sec:model_problem} is dedicated to the detailed study of model problems in multiobjective bilevel optimization
with a particular emphasis on value-function-type reformulations and their relations to each other.
Based on a discussion of closedness properties of efficiency-type and associated frontier mappings in
\cref{sec:param_vec_opt}, some existence results are stated in \cref{sec:existence}. 
Furthermore, an intermediate way between efficiency and weak efficiency on how to interpret \eqref{eq:parametric_problem}
is suggested in \cref{sec:intermediate_approach}.
In \cref{sec:coderivative_calculus}, we are concerned with the derivation of upper coderivative estimates for
frontier mappings. Based on a classical approach from the literature, we recover some known results in the novel setting
of weak efficiency in \cref{sec:coderivative_domination}. Relying on a scalarization approach, some new findings are
obtained in \cref{sec:coderivative_scalarization}.
The paper closes with some concluding remarks in \cref{sec:conclusions}.

\section{Preliminaries}\label{sec:preliminaries}

In this section, we comment on the notation used in this paper.
Particularly, some foundations of set-valued and variational analysis are provided.
Furthermore, we present some necessary basics of multiobjective optimization.

\subsection{Basic notation}

The notation in this paper is fairly standard and mainly follows
\cite{Mordukhovich2006,RockafellarWets2009}.

Let $\R_+$ and $\R_-$ represent the nonnegative and nonpositive real numbers, respectively.
Throughout the paper, we equip $\R^n$, where $n\in\N$ is a positive integer, 
with the Euclidean inner product
$\langle \cdot,\cdot\rangle\colon\R^n\times\R^n\to\R$
and the associated Euclidean norm
$\norm{\cdot}\colon\R^n\to\R_+$.
We use $\mathbb S_1(0)$ to represent the unit sphere around the origin in $\R^n$ where the
dimension of the underlying space will be always clear from the context.
For some set $\Omega\subset\R^n$, $\cl\Omega$, $\intr\Omega$, $\conv\Omega$, and $\cone\Omega$
denote the closure, the interior, the convex hull, and the convex conic hull of $\Omega$,
respectively.

Whenever $\upsilon\colon\R^n\to\R^m$ is a mapping which is differentiable at $\bar x\in\R^n$,
then $\upsilon'(\bar x)\in\R^{m\times n}$ denotes the Jacobian of $\upsilon$ at $\bar x$.
Partial derivatives are represented in analogous fashion.

\subsubsection{Convexity properties of vector-valued functions}

In this paper, we are partially concerned with convexity properties of set-valued mappings
with respect to (w.r.t.) a given so-called ordering cone, see e.g.\ \cite{Jahn2011}.
For a nonempty, closed, convex, pointed cone $\mathcal K\subset\R^p$, a mapping
$F\colon\R^m\to\R^p$ is said to be $\mathcal K$-convex if
\[
	\forall y,y'\in\R^m\,\forall\alpha\in[0,1]\colon\quad
	\alpha F(y)+(1-\alpha)F(y')-F(\alpha y+(1-\alpha)y')\in \mathcal K.
\]
We note that whenever $F$ is $\mathcal K$-convex and $Y\subset\R^m$ is convex,
then $F(Y)+\mathcal K$ is convex,
and the function $y\mapsto\langle \lambda,F(y)\rangle$ is convex for each $\lambda\in \mathcal K^*$.
Here, we made use of 
\[
	\mathcal K^*:=\{\lambda\,|\,\forall z\in \mathcal K\colon\,\langle\lambda,z\rangle\geq 0\}
\]
to denote the dual cone of $\mathcal K$ which is a closed, convex cone as well.
Let us mention that a componentwise convex mapping $F\colon\R^m\to\R^p$ is $\R^p_+$-convex.

The mapping $F$ is said to be $\mathcal K$-strictly-quasiconvex if 
\begin{align*}
	\forall y,y' \in \R^m\,&\forall\alpha\in(0,1)\colon\\
 	&F(y')-F(y) \in \mathcal K,\,y\neq y' \quad\Longrightarrow\quad F(y') -   F(\alpha y + (1 - \alpha) y') \in \intr \mathcal K.  
\end{align*}
Exemplary, in the case where $\mathcal K := \R^p_+$, $F$ is $\R^p_+$-strictly-quasiconvex 
if $F_1,\ldots,F_p$ are so-called semi-strictly quasiconvex in the sense that
\[
	\forall y,y'\in\R^m\,\forall \alpha\in(0,1)\colon\quad
	F_i(y')\geq F_i(y),\,y\neq y'
	\quad
	\Longrightarrow
	\quad
	F_i(y')> F_i(\alpha y+(1-\alpha)y')
\]
has to hold
for all $i=1,\ldots,p$. The latter is clearly inherent whenever the functions
$F_1,\ldots,F_p$ are strictly convex.

\subsubsection{Set-valued mappings}

In this paper, we are concerned with so-called set-valued mappings which assign to each
element of the pre-image space a (potentially empty) subset of the image space.
Such objects will be denoted by $\Upsilon\colon\R^n\tto\R^m$.
Throughout the paper, we use 
\[
	\gph\Upsilon:=\{(x,y)\,|\,y\in\Upsilon(x)\},
	\qquad
	\dom\Upsilon:=\{x\,|\,\Upsilon(x)\neq\emptyset\}
\] 
to represent the graph and the domain of $\Upsilon$, respectively.
For a set $S\subset\R^m$, we define the set-valued mappings 
$\Upsilon+S,\Upsilon-S\colon\R^n\tto\R^m$ by means of
\[
	\forall x\in\R^n\colon\quad
	(\Upsilon+S)(x):=\Upsilon(x)+S,\qquad
	(\Upsilon-S)(x):=\Upsilon(x)-S
\]
for brevity of notation.
A set-valued mapping is referred to as polyhedral whenever its graph can
be represented as the union of finitely many convex polyhedral sets.

For fixed $\bar x\in\dom\Upsilon$, we call $\Upsilon$ closed at $\bar x$
whenever for each sequence $((x^k,y^k))_{k\in\N}\subset\gph\Upsilon$ and each
$\bar y\in\R^m$ such that $x^k\to\bar x$ and $y^k\to\bar y$, 
we also have $\bar y\in\Upsilon(\bar x)$. Note that this property is inherent whenever
$\gph\Upsilon$ is closed. Additionally, whenever $\Upsilon$ is closed at $\bar x$,
$\Upsilon(\bar x)$ is a closed set.
Observe that closedness of $\gph\Upsilon$ does not necessarily give closedness of $\dom\Upsilon$.
However, closedness of $\dom\Upsilon$ and closedness of $\Upsilon$ at each point
from $\dom\Upsilon$ give closedness of $\gph\Upsilon$.
We refer to $\Upsilon$ as locally bounded at $\bar x$ if there exist
a bounded set $B\subset\R^m$ and some neighborhood $U\subset\R^n$ of $\bar x$ 
such that $\Upsilon(x)\subset B$
for all $x\in U$.
Furthermore, $\Upsilon$ is called lower semicontinuous at $\bar x$
(in Berge's sense) w.r.t.\ $\Omega\subset\R^n$ whenever for each open set $V\subset\R^m$
such that $\Upsilon(\bar x)\cap V\neq\emptyset$, there is neighborhood $U\subset\R^n$ of $\bar x$ such that
$\Upsilon(x)\cap V\neq\emptyset$ for all $x\in\Omega\cap U$.
For some fixed pair $(\bar x,\bar y)\in\gph\Upsilon$, we say that $\Upsilon$
is inner semicontinuous at $(\bar x,\bar y)$ w.r.t.\ $\Omega\subset\R^n$ whenever for
each neighborhood $V\subset\R^m$ of $\bar y$, there exists a neighborhood $U\subset\R^n$ of $\bar x$ such that
$\Upsilon(x)\cap V\neq\emptyset$ for all $x\in\Omega\cap U$.
If $\Omega:=\R^n$ can be chosen, $\Upsilon$ is called inner semicontinuous at $(\bar x,\bar y)$
for brevity.
Note that $\Upsilon$ is lower semicontinuous at $\bar x$ w.r.t.\ $\Omega$
if and only if it is inner semicontinuous at $(\bar x,y)$ w.r.t.\ $\Omega$
for each $y\in\Upsilon(\bar x)$.
Furthermore, observe that $\Upsilon$ is inner semicontinuous at $(\bar x,\bar y)$ w.r.t.\ $\Omega$
if and only if for each sequence $(x^k)_{k\in\N}\subset\Omega$ such that $x^k\to\bar x$,
there exists a sequence $(y^k)_{k\in\N}\subset\R^m$ satisfying $y^k\to\bar y$ and
$y^k\in\Upsilon(x^k)$ for all large enough $k\in\N$.
Let us also mention that $\Upsilon$ is referred to as inner semicompact at $\bar x$ w.r.t.\ $\Omega$
if for each sequence $(x^k)_{k\in\N}\subset\Omega$ with $x^k\to\bar x$, there are a subsequence
$(x^{k_\ell})_{\ell\in\N}$ and a convergent sequence $(y^\ell)_{\ell\in\N}$ such that
$y^\ell\in\Upsilon(x^{k_\ell})$ holds for all $\ell\in\N$.
If $\Omega:=\R^n$ can be chosen, $\Upsilon$ is called inner semicompact at $\bar x$.
Clearly, if $\Upsilon$ is locally bounded at $\bar x$, then it is inner semicompact at $\bar x$
w.r.t.\ each subset of $\R^n$. Furthermore, whenever there exists $\bar y\in\Upsilon(\bar x)$ such that
$\Upsilon$ is inner semicontinuous at $(\bar x,\bar y)$ w.r.t.\ $\Omega$, then $\Upsilon$ is
inner semicompact at $\bar x$ w.r.t.\ $\Omega$.
In the definitions of lower and inner semicontinuity as well as inner semicompactness, 
special emphasis is laid on settings where $\Omega\in\{\R^n,\dom\Upsilon\}$.

In the subsequently stated lemma, we characterize the closedness of the graph of set-valued mappings of special type.

\begin{lemma}\label{lem:closed_graph_for_special_map}
	For a function $\upsilon\colon\R^n\to\overline\R$, we consider the set-valued mapping
	$\Upsilon\colon\R^n\tto\R$ given by 
	\[
		\forall x\in\R^n\colon\quad
		\Upsilon(x)
		:=
		\{\upsilon(x)\}-\R_+
	\]
	with the conventions $\{\infty\}-\R_+:=\R$ and $\{-\infty\}-\R_+:=\emptyset$.
	Then $\gph\Upsilon$ is closed if and only if $\upsilon$ is upper semicontinuous.
\end{lemma}
\begin{proof}
		We observe that $\gph\Upsilon$ equals 
		$\operatorname{hypo}\upsilon:=\{(x,\alpha)\in\R^n\times\R\,|\,\alpha\leq \upsilon(x)\}$,
		the so-called hypograph of $\upsilon$, 
		and $\operatorname{hypo}\upsilon$ is closed if and only if $\upsilon$ is an
		upper semicontinuous function, see e.g.\ \cite[Theorem~1.6]{RockafellarWets2009}.
\end{proof}

For two given set-valued mappings $\Upsilon_1\colon\R^n\tto\R^m$ and $\Upsilon_2\colon\R^m\tto\R^p$,
their composition $\Upsilon_2\circ\Upsilon_1\colon\R^n\tto\R^p$ is defined by
\[
	\forall x\in\R^n\colon\quad
	(\Upsilon_2\circ\Upsilon_1)(x)
	:=
	\bigcup\limits_{y\in\Upsilon_1(x)}\Upsilon_2(y).
\]
The associated so-called intermediate mapping $\Xi\colon\R^n\times\R^p\tto\R^m$ is given by
\[
	\forall x\in\R^n\,\forall z\in\R^p\colon\quad
	\Xi(x,z):=\{y\in\Upsilon_1(x)\,|\,z\in\Upsilon_2(y)\}
\]
and plays an important role for the analysis of the underlying composition,
see e.g.\ \cite[Section~5.3]{BenkoMehlitz2022} or \cite[Section~3.1.2]{Mordukhovich2006}.

\subsubsection{Variational analysis and generalized differentiation}\label{sec:VA}

The following definitions are taken from \cite{Mordukhovich2006,RockafellarWets2009}.

For some set $\Omega\subset\R^n$, which is closed locally around some point $\bar x\in\Omega$, we call
\[
	N_\Omega(\bar x)
	:=
	\left\{
		\eta\,\middle|\,
			\begin{aligned}
				&\exists(x^k)_{k\in\N},(y^k)_{k\in\N}\subset\R^n\,\exists(\alpha_k)_{k\in\N}\subset(0,\infty)\colon\\
				&\quad x^k\to\bar x,\,\alpha_k(x^k-y^k)\to\eta,\,y^k\in\Pi_\Omega(x^k)\,\forall k\in\N
			\end{aligned}
	\right\}
\]
the limiting (or Mordukhovich) normal cone to $\Omega$ at $\bar x$.
Above, $\Pi_\Omega\colon\R^n\tto\R^n$ denotes the projector of $\Omega$ which is given by
\[
	\forall x\in\R^n\colon\quad\Pi_\Omega(x):=\argmin\limits_y\{\norm{y-x}\,|\,y\in\Omega\}.
\]
We note that $N_\Omega(\bar x)$ is a closed cone which reduces to the classical normal cone
of convex analysis as soon as $\Omega$ is a convex set.

For a set-valued mapping $\Upsilon\colon\R^n\tto\R^m$ and some point $(\bar x,\bar y)\in\gph\Upsilon$
where $\gph\Upsilon$ is locally closed, the set-valued mapping 
$D^*\Upsilon(\bar x,\bar y)\colon\R^m\tto\R^n$ given by
\[
	\forall y^*\in\R^m\colon\quad
	D^*\Upsilon(\bar x,\bar y)(y^*):=\{x^*\,|\,(x^*,-y^*)\in N_{\gph\Upsilon}(\bar x,\bar y)\}
\]
is referred to as the (limiting or Mordukhovich) coderivative of $\Upsilon$ at $(\bar x,\bar y)$.
In the case of a single-valued continuous mapping $\upsilon\colon\R^n\to\R^m$, 
which can naturally be interpreted as a singleton-valued set-valued mapping, 
and $\bar x\in\R^n$, we make use
of $D^*\upsilon(\bar x):=D^*\Upsilon(\bar x,\upsilon(\bar x))$ for brevity of notation.
If the mapping $\upsilon$ is continuously differentiable at $\bar x$, then
$D^*\upsilon(\bar x)(y^*)=\{\upsilon'(\bar x)^\top y^*\}$ is valid for all $y^*\in\R^m$.

\subsection{Foundations of multiobjective optimization}\label{sec:vector_optimization}

For a continuous vector function $F\colon\R^m\to\R^p$, a nonempty, closed set $S\subset\R^m$, and some closed, convex, pointed
cone $\mathcal K\subset\R^p$ with nonempty interior, we investigate the rather general multiobjective optimization problem
\begin{equation}\label{eq:MOP}\tag{MOP}
	\min{}_{\mathcal{K}}\{F(y)\,|\,y\in S\}.
\end{equation}
Recall that some feasible point $\bar y\in S$ of \eqref{eq:MOP} is called efficient (weakly efficient) whenever
\[
	(\{F(\bar y)\}-\mathcal K\setminus\{0\})\cap F(S)=\emptyset
	\qquad
	\bigl((\{F(\bar y)\}-\intr \mathcal K)\cap F(S)=\emptyset\bigr)
\]
holds. The set of all efficient (weakly efficient) points of \eqref{eq:MOP} will be denoted by
$\eff{F}{S}{\mathcal K}$ ($\weff{F}{S}{\mathcal K}$). 
Further, $\nd{F}{S}{\mathcal K}:=F(\eff{F}{S}{\mathcal K})$ 
($\wnd{F}{S}{\mathcal K}:=F(\weff{F}{S}{\mathcal K})$) is referred to as the nondominated
(weakly nondominated) set of \eqref{eq:MOP}. Let us mention that the computation of efficient (weakly efficient) points
makes knowledge of $F$, $S$, and $\mathcal K$ necessary, while, in principle, nondominated (weakly nondominated) points can
be characterized from knowledge of $F(S)$ and $\mathcal K$ only.
Clearly, the inclusions $\eff{F}{S}{\mathcal K}\subset\weff{F}{S}{\mathcal K}$ 
and $\nd{F}{S}{\mathcal K}\subset\wnd{F}{S}{\mathcal K}$ hold by definition.
Let us also point out that continuity of $F$ and the openness of $\intr \mathcal K$ 
give closedness of the sets $\weff{F}{S}{\mathcal K}$ and $\wnd{F}{S}{\mathcal K}$.
In contrast, the sets $\eff{F}{S}{\mathcal K}$ and $\nd{F}{S}{\mathcal K}$ do not need to be closed, 
see \cref{ex:domination_property_without_closedness} below.
Let us also note that the closedness of $\weff{F}{S}{\mathcal K}$ and $\wnd{F}{S}{\mathcal K}$ guarantees
$\cl\eff{F}{S}{\mathcal K}\subset\weff{F}{S}{\mathcal K}$ 
and $\cl\nd{F}{S}{\mathcal K}\subset\wnd{F}{S}{\mathcal K}$, and 
these inclusions can be strict.
In the case where $p:=1$, there is no difference between efficiency and weak efficiency as well as
nondominance and weak nondominance, respectively, since the only cones satisfying our
requirements would be $\R_+$ and $\R_-$ in this case.

In this paper, we are concerned with so-called domination properties of \eqref{eq:MOP},
see e.g.\ \cite{Luc1984,Henig9867,LiuNgYang2009}.

\begin{definition}\label{def:domination}
	We say that \eqref{eq:MOP} possesses the domination property (the weak domination property)
	if the following inclusion is valid:
	\[
		F(S)\subset\nd{F}{S}{\mathcal K}+\mathcal K
		\qquad
		\Big(F(S)\subset\wnd{F}{S}{\mathcal K}+\mathcal K\Big).
	\]
\end{definition}

It is clear from the definition that validity of the domination property also gives validity of the
weak domination property.
Observe that, by convexity of $\mathcal K$, \eqref{eq:MOP} possesses the domination property (the weak domination property)
if and only if
\[
	F(S)+\mathcal K=\nd{F}{S}{\mathcal K}+\mathcal K
	\qquad
	\Big(F(S)+\mathcal K=\wnd{F}{S}{\mathcal K}+\mathcal K\Big).
\]
Thus, the discussion at the beginning of \cite[Section~4]{LafhimZemkoho2022} is completely superfluous.
Sufficient criteria for the domination and weak domination property can be found in \cite{Henig9867,LiuNgYang2009}. 
Exemplary, let us mention that \eqref{eq:MOP} possesses the domination property whenever $F$ is continuous while
$S$ is compact. This follows by combining \cite[Theorem~2]{Henig9867} and \cite[Theorem~3.38]{Jahn2011}.
Furthermore, we also want to point the reader's attention to the fact that even in the presence of the domination property,
$\eff{F}{S}{\mathcal K}$ and $\nd{F}{S}{\mathcal K}$ might be non-closed.

\begin{example}\label{ex:domination_property_without_closedness}
	Let $F\colon\R^2\to\R^2$ be the identity, let $S\subset\R^2$ be given by
	\[
		S:=\{y\in\R^2_+\,|\,y_1^2+y_2^2\geq 1\}\cup([-1,0]\times[1,\infty)),
	\]
	and let $\mathcal K:=\R^2_+$. 
	Then we find
	\[
		\eff{F}{S}{\mathcal K}
		=
		\nd{F}{S}{\mathcal K}
		=
		\{(\cos t,\sin t)\in\R^2\,|\,t\in[0,\pi/2)\}\cup\{(-1,1)\}.
	\]
	On the one hand, this can be used to see that the associated problem \eqref{eq:MOP} possesses the
	domination property.
	On the other hand, $\eff{F}{S}{\mathcal K}$ and $\nd{F}{S}{\mathcal K}$ are not closed.
\end{example}

It follows from \cite[Example~3.4.2]{SawaragiNakayamaTanino1985} that the efficient set $\eff{F}{S}{\mathcal K}$ 
might not be closed even if the set $F(S)$ is compact and convex. 
However, for instance in \cite[Lemma 7.1]{Plastria2020}, it is shown that if $\mathcal K := \R^p_+$ and all functions $F_1,\ldots,F_p$ 
are strictly quasiconvex while $S$ is convex, then $\eff{F}{S}{\mathcal K}$ coincides with $\weff{F}{S}{\mathcal K}$ 
and so $\eff{F}{S}{\mathcal K}$ is closed since $\weff{F}{S}{\mathcal K}$ is closed. 
This result can also be extended to more general ordering cones. 

\begin{lemma}\label{lem:closedness_of_efficient_set}
Let $F$ be $\mathcal K$-strictly-quasiconvex, and let $S$ be convex. 
Then $\eff{F}{S}{\mathcal K} = \weff{F}{S}{\mathcal K}$ and $\nd{F}{S}{\mathcal K}=\wnd{F}{S}{\mathcal K}$.
\end{lemma}
\begin{proof}
	We only show equality of the efficient and weakly efficient set.
	This automatically gives equivalence of the nondominated and weakly nondominated set.
	Since $\eff{F}{S}{\mathcal K} \subset \weff{F}{S}{\mathcal K}$ holds in general, we only need to prove the converse inclusion. 
	Fix $\bar{y} \in \weff{F}{S}{\mathcal K}$ and suppose that $\bar{y} \notin \eff{F}{S}{\mathcal K}$. 
	Then we find some $\tilde{y} \in S \setminus \{ \bar y \}$ satisfying 
	$F(\tilde{y}) \in \{F(\bar{y})\} - \mathcal K \setminus \{ 0 \}$. 
	For all $\alpha \in (0,1)$, $\alpha \tilde{y} + (1 - \alpha) \bar{y} \in S$ follows by convexity of $S$,
	and since $F$ is $\mathcal K$-strictly-quasiconvex, we obtain
	\begin{align*}
		F(\alpha \tilde{y} + (1 - \alpha) \bar{y}) \in (\{F(\bar{y})\} - \intr \mathcal K) \cap F(S). 
	\end{align*}
	This is a contradiction to $\bar{y} \in \weff{F}{S}{\mathcal K}$. 
\end{proof}

\section{Value function models in multiobjective bilevel optimization}\label{sec:model_problem}

In order to make \eqref{eq:BPP} explicit, one has to specify the precise meaning behind the
optimization goals at the (upper- and) lower-level stage.
Therefore, let us define set-valued mappings 
$\Psi,\Psi_\textup{w}\colon\R^n\tto\R^m$
and 
$\Phi,\Phi_\textup{w}\colon\R^n\tto\R^q$ by means of
\begin{equation*}
	\begin{aligned}
		\Psi(x)&:=\eff{f(x,\cdot)}{\Gamma(x)}{\mathcal C},&
		\qquad
		\Phi(x)&:=\nd{f(x,\cdot)}{\Gamma(x)}{\mathcal C},&\\
		\Psi_\textup{w}(x)&:=\weff{f(x,\cdot)}{\Gamma(x)}{\mathcal C},&
		\qquad
		\Phi_\textup{w}(x)&:=\wnd{f(x,\cdot)}{\Gamma(x)}{\mathcal C}&
	\end{aligned}
\end{equation*}
for each $x\in\R^n$. 
Note that the identities
\[
	\dom\Psi=\dom\Phi,\qquad
	\dom\Psi_\textup{w}=\dom\Phi_\textup{w}
\]
hold by definition of these mappings.
We will refer to $\Psi$ and $\Psi_\textup{w}$ as the efficiency and weak efficiency mapping
associated with \eqref{eq:parametric_problem}, respectively.
Similarly, $\Phi$ and $\Phi_\textup{w}$ are called the frontier and weak frontier mapping 
associated with \eqref{eq:parametric_problem}.
For brevity of notation, we further introduce the set-valued mapping 
$\Sigma\colon\R^n\tto\R^q$ given by
\begin{equation}\label{eq:surrogate_maps_coderivative_calculus}
	\forall x\in\R^n\colon\quad
	\Sigma(x):=f(x,\Gamma(x)).
\end{equation}

In the course of this note, we are mainly concerned with the optimization problems
\begin{subequations}
	\begin{align}
	\label{eq:BPP_eff}\tag{$\mathcal E$-BPP}
	&\min\limits_{x,y}{}_{\mathcal{K}}\{F(x,y)\,|\,x\in X,\,y\in\Psi(x)\},
	\\
	\label{eq:BPP_weff}\tag{$\mathcal E_\textup{w}$-BPP}
	&\min\limits_{x,y}{}_{\mathcal{K}}\{F(x,y)\,|\,x\in X,\,y\in\Psi_\textup{w}(x)\}
	\end{align}
\end{subequations}
as well as their associated respective value function reformulations
\begin{subequations}
	\begin{align}
	\label{eq:VFR_eff}\tag{$\mathcal E$-VFR}
	&\min\limits_{x,y}{}_{\mathcal{K}}\{F(x,y)\,|\,x\in X,\,f(x,y)\in\Phi(x),\,y\in\Gamma(x)\},
	\\
	\label{eq:VFR_weff}\tag{$\mathcal E_{\textup{w}}$-VFR}
	&\min\limits_{x,y}{}_{\mathcal{K}}\{F(x,y)\,|\,x\in X,\,f(x,y)\in\Phi_\textup{w}(x),\,y\in\Gamma(x)\}
	\end{align}
\end{subequations}
which are apparently equivalent. 

\begin{proposition}\label{prop:value_function_char_of_eff_set}
	\leavevmode
	\begin{enumerate}
		\item For each $x\in\R^n$, $\Psi(x)=\{y\in\Gamma(x)\,|\,f(x,y)\in\Phi(x)\}$ holds true.
		\item For each $x\in\R^n$, $\Psi_\textup{w}(x)=\{y\in\Gamma(x)\,|\,f(x,y)\in\Phi_\textup{w}(x)\}$ holds true.
	\end{enumerate}
\end{proposition}
\begin{proof}
	Let us start with the proof of the first assertion.
	Indeed, if $y\in\Psi(x)$ holds for some $x\in\R^n$, then $y\in\Gamma(x)$ and $f(x,y)\in f(x,\Psi(x))=\Phi(x)$
	follow trivially, i.e., the inclusion $\subset$ is obvious.
	For the proof of the converse inclusion $\supset$, pick $y\in\Gamma(x)$ such that $f(x,y)\in\Phi(x)$ for some $x\in\R^n$.
	By definition of $\Phi$, we find $y'\in\Psi(x)$ such that $f(x,y)=f(x,y')$.
	Since $(\{f(x,y')\}-\mathcal C\setminus\{0\})\cap \Sigma(x)=\emptyset$, we also have
	$(\{f(x,y)\}-\mathcal C\setminus\{0\})\cap \Sigma(x)=\emptyset$, and this gives $y\in\Psi(x)$.
	
	The proof of the second assertion is completely analogous.
\end{proof}

\begin{corollary}\label{cor:VFR_equivelant}
	\leavevmode
	\begin{enumerate}
		\item Problems \eqref{eq:BPP_eff} and \eqref{eq:VFR_eff} possess the same feasible set.
		\item Problems \eqref{eq:BPP_weff} and \eqref{eq:VFR_weff} possess the same feasible set.
	\end{enumerate}
\end{corollary}

Let us note that whenever $q:=1$ and $\mathcal C:=\R_+$ hold, 
then both problems \eqref{eq:VFR_eff} and \eqref{eq:VFR_weff} are equivalent to
\begin{equation}\label{eq:VFR_nonstandard}
	\min\limits_{x,y}{}_{\mathcal{K}}\{F(x,y)\,|\,x\in X,\,f(x,y)=\varphi(x),\,y\in\Gamma(x)\}
\end{equation}
which, up to the relation sign in the constraint involving $\varphi$, is the same as
the classical value function reformulation \eqref{eq:VFR_classic}.
In scalar bilevel programming, however, it has turned out to be beneficial to work
with the constraint $f(x,y)\leq\varphi(x)$ as this opens the door to algorithmic
applications and induces a sign condition for the Lagrange multiplier associated 
with this constraint in the context of nonsmooth optimality conditions.
Noting that, for each $x\in\R^n$ and $y\in\Gamma(x)$, $f(x,y)<\varphi(x)$ is not possible,
we can reformulate \eqref{eq:VFR_classic} and \eqref{eq:VFR_nonstandard} as
\begin{equation}\label{eq:VFR_classic_with_ordering_cone}
	\min\limits_{x,y}{}_{\mathcal{K}}\{F(x,y)\,|\,x\in X,\,f(x,y)\in\{\varphi(x)\}-\R_+,\,y\in\Gamma(x)\},
\end{equation}
with the conventions $\{\infty\}-\R_+:=\R$ and $\{-\infty\}-\R_+:=\emptyset$.
Observing that $\R_+$ is the ordering cone associated with the lower-level problem 
in the scalar situation, this motivates to investigate
\begin{subequations}
	\begin{align}
	\label{eq:VFR_eff_C}
	&\min\limits_{x,y}{}_{\mathcal{K}}\{F(x,y)\,|\,x\in X,\,f(x,y)\in\Phi(x)-\mathcal C,\,y\in\Gamma(x)\},
	\\
	\label{eq:VFR_weff_C}
	&\min\limits_{x,y}{}_{\mathcal{K}}\{F(x,y)\,|\,x\in X,\,f(x,y)\in\Phi_\textup{w}(x)-\mathcal C,\,y\in\Gamma(x)\}
	\end{align}
\end{subequations}
as potential surrogates of 
\eqref{eq:VFR_eff} and \eqref{eq:VFR_weff}, respectively. 
As the following results show, the addition of the lower-level ordering cone 
in \eqref{eq:VFR_eff_C} and \eqref{eq:VFR_weff_C} does not change the problem.

\begin{proposition}\label{prop:value_function_char_extended_of_eff_set}
	\leavevmode
	\begin{enumerate}
		\item For each $x\in\R^n$, $\Psi(x)=\{y\in\Gamma(x)\,|\,f(x,y)\in\Phi(x)-\mathcal C\}$ holds true.
		\item For each $x\in\R^n$, $\Psi_\textup{w}(x)=\{y\in\Gamma(x)\,|\,f(x,y)\in\Phi_\textup{w}(x)-\mathcal C\}$ holds true.
	\end{enumerate}
\end{proposition}
\begin{proof}
	Let us start with the proof of the first assertion.
	Clearly, as $0\in \mathcal C$, the inclusion $\subset$ is clear from \cref{prop:value_function_char_of_eff_set}.
	Conversely, fix $x\in\R^n$ and $y\in\Gamma(x)$ such that $f(x,y)\in\Phi(x)-\mathcal C$.
	Then there exist $y'\in\Psi(x)$ and $c\in \mathcal C$ such that $f(x,y)=f(x,y')-c$.
	By definition of efficiency, this gives $c=0$, i.e., $f(x,y)=f(x,y')\in\Phi(x)$, and, hence,
	$y\in\Psi(x)$ by \cref{prop:value_function_char_of_eff_set}.
	
	In order to prove the second assertion, we observe that the inclusion $\subset$ follows again
	from $0\in \mathcal C$ and \cref{prop:value_function_char_of_eff_set}.
	Thus, pick an arbitrary $x\in\R^n$ and $y\in\Gamma(x)$ such that $f(x,y)\in\Phi_\textup{w}(x)-\mathcal C$.
	Then there exist $y'\in\Psi_\textup{w}(x)$ and $c\in \mathcal C$ such that $f(x,y)=f(x,y')-c$.
	Suppose that $y\notin\Psi_\textup{w}(x)$. Then there are $\tilde y\in\Gamma(x)$ and $\tilde c\in\intr \mathcal C$
	such that $f(x,\tilde y)=f(x,y)-\tilde c$, and this gives 
	$f(x,\tilde y)=f(x,y')-(c+\tilde c)$. As we have $c+\tilde c\in \intr \mathcal C$, $y'\notin\Psi_\textup{w}(x)$
	follows which is a contradiction. Hence, we find $y\in\Psi_\textup{w}(x)$ which yields
	$f(x,y)\in\Phi_\textup{w}(x)$, i.e., $y\in\Psi_\textup{w}(x)$ due to \cref{prop:value_function_char_of_eff_set}.
\end{proof}

\begin{corollary}\label{cor:adding_an_ordering_cone_in_VFR}
	\leavevmode
	\begin{enumerate}
		\item Problems \eqref{eq:VFR_eff} and \eqref{eq:VFR_eff_C} possess the same feasible set.
		\item Problems \eqref{eq:VFR_weff} and \eqref{eq:VFR_weff_C} possess the same feasible set.
	\end{enumerate}
\end{corollary}

In order to guarantee the existence of efficient or weakly efficient points for the above
model problems \eqref{eq:BPP_eff}, \eqref{eq:BPP_weff}, \eqref{eq:VFR_eff}, and \eqref{eq:VFR_weff}, 
and to make their algorithmic treatment reasonable, it is important to
guarantee that their respective feasible sets are nonempty and closed.
Recalling that the mapping $f$ is continuous, and the sets $X$ and $\gph\Gamma$ are assumed to
be closed, this reduces to the closedness of $\gph\Psi$, $\gph\Psi_\textup{w}$,
$\gph\Phi$, and $\gph\Phi_\textup{w}$, respectively.
This is addressed in \cref{sec:param_vec_opt,sec:existence} below.
As we will see, $\gph\Psi$ and $\gph\Phi$ are not likely to be closed in several popular settings.
In \cref{sec:intermediate_approach}, we suggest an intermediate approach that bridges the
model problems \eqref{eq:BPP_eff} and \eqref{eq:BPP_weff} 
(and, similarly, \eqref{eq:VFR_eff} and \eqref{eq:VFR_weff}) and comes along with a
guaranteed closedness of the feasible set.

\subsection{Closedness properties of efficiency and frontier mappings}\label{sec:param_vec_opt}

In the following lemma, we show that $\Psi_\textup{w}$
as well as $\Phi_\textup{w}$ are indeed closed under mild assumptions,
see e.g.\ \cite[Corollary~3.5.6]{GoepfertRiahiTammerZalinescu2003} for a related result.

\begin{lemma}\label{thm:closedness_of_graphs}
	Fix $\bar x\in\dom\Psi_\textup{w}$ and let $\Gamma$ be lower semicontinuous 
	at $\bar x$ w.r.t.\ its domain. Then $\Psi_\textup{w}$ is closed at $\bar x$.
	If, additionally, $\Gamma$ is locally bounded at $\bar x$, $\Phi_\textup{w}$
	is closed at $\bar x$.
\end{lemma}
\begin{proof}
	Let $(x^k)_{k\in\N}\subset\R^n$ and $(y^k)_{k\in\N}\subset\R^m$ as well as 
	$\bar y\in\R^m$ be chosen such that $x^k\to\bar x$, $y^k\to\bar y$, and
	$y^k\in\Psi_\textup{w}(x^k)$ for all $k\in\N$. By definition of $\Psi_\textup{w}$,
	we find $y^k\in\Gamma(x^k)$ and
	\begin{equation}\label{eq:parametric_weak_efficiency}
		(\{f(x^k,y^k)\}-\intr \mathcal C)\cap \Sigma(x^k)=\emptyset
	\end{equation}
	for all $k\in\N$. 
	Since $\Gamma$ is closed at $\bar{x}$, we obtain $\bar y\in\Gamma(\bar x)$.
	Suppose that $\bar y\notin\Psi_\textup{w}(\bar x)$. 
	Then there are $\tilde y\in\Gamma(\bar x)$ and $\bar z\in\intr \mathcal C$ such that
	$f(\bar x,\bar y)-\bar z=f(\bar x,\tilde y)$.
	Noting that $(x^k)_{k\in\N}\subset\dom\Psi_\textup{w}\subset\dom\Gamma$,
	the assumed lower semicontinuity of $\Gamma$ yields the existence of
	a sequence $(\tilde y^k)_{k\in\N}\subset\R^m$ such that $\tilde y^k\to\tilde y$ and
	$\tilde y^k\in\Gamma(x^k)$ for all $k\in\N$.
	We set
	\[
		z^k:=f(x^k,y^k)-f(\bar x,\bar y)-f(x^k,\tilde y^k)+f(\bar x,\tilde y)+\bar z
	\]
	for each $k\in\N$. Continuity of $f$ gives us $z^k\to\bar z$, 
	and due to $\bar z\in\intr \mathcal C$, $z^k\in\intr \mathcal C$ holds for all large enough $k\in\N$.
	Some rearrangements yield
	\[
		f(x^k,y^k)-z^k
		=
		f(x^k,\tilde y^k)+\bigl(f(\bar x,\bar y)-\bar z-f(\bar x,\tilde y)\bigr)
		=
		f(x^k,\tilde y^k),
	\]
	i.e.,
	\[
		(\{f(x^k,y^k)\}-\intr \mathcal C)\cap \Sigma(x^k)\neq\emptyset
	\]
	for large enough $k\in\N$, contradicting \eqref{eq:parametric_weak_efficiency}.
	Hence, $\bar y\in\Psi_\textup{w}(\bar x)$ follows, 
	i.e., $\Psi_\textup{w}$ is closed at $\bar x$.
	
	In order to verify the closedness of $\Phi_\textup{w}$ at $\bar x$, we pick
	sequences $(x^k)_{k\in\N}\subset\R^n$ and $(z^k)_{k\in\N}\subset\R^p$ such that
	$x^k\to\bar x$, $z^k\in\Phi_\textup{w}(x^k)$ for all $k\in\N$, and $z^k\to\bar z$
	for some $\bar z\in\R^p$. This gives the existence of $(y^k)_{k\in\N}\subset\R^m$
	such that $y^k\in\Psi_\textup{w}(x^k)$ and $f(x^k,y^k)=z^k$ for all $k\in\N$.
	Local boundedness of $\Gamma$ at $\bar x$ guarantees $y^k\to\bar y$ along a subsequence
	(without relabeling) for some $\bar y\in \R^m$.
	The first part of the proof ensures $\bar y\in\Psi_\textup{w}(\bar x)$, and continuity of
	$f$ yields $\bar z=f(\bar x,\bar y)$, i.e., $\bar z\in\Phi_\textup{w}(\bar x)$,
	showing the closedness of $\Phi_\textup{w}$ at $\bar x$.
\end{proof}

As a corollary, we can state conditions which guarantee 
the closedness of $\gph\Psi_\textup{w}$ and $\gph\Phi_\textup{w}$.

\begin{corollary}\label{cor:closedness_of_graphs}
	Let $\dom\Gamma$ be closed, and let $\Gamma$ be locally bounded and lower semicontinuous
	w.r.t.\ its domain at each point from $\dom\Gamma$. Then $\gph\Psi_{\textup{w}}$ and $\gph\Phi_\textup{w}$
	are closed.
\end{corollary}
\begin{proof}
	The closedness of $\gph\Gamma$ and the assumed local boundedness guarantee that,
	for each $\bar x\in\dom\Gamma$, $\Gamma(\bar x)$ is nonempty and compact.
	Thus, \cite[Theorem~6.3]{Jahn2011} shows $\dom\Psi_\textup{w}=\dom\Gamma$,
	i.e., $\dom\Psi_\textup{w}$ and $\dom\Phi_\textup{w}$ are closed.
	Furthermore, \cref{thm:closedness_of_graphs} shows that $\Psi_\textup{w}$ and $\Phi_\textup{w}$
	are closed at each point of their domain, respectively.
	Taking these facts together, the assertion follows.
\end{proof}

Let us relate these findings to the classical results from \cite{BankGuddatKlatteKummerTammer1983}
which apply to
the scalar situation $q:=1$ and $\mathcal C:=\R_+$. 
Then the closedness of $\Psi_\textup{w}$ corresponds to
the closedness of the classical solution mapping while the closedness of $\Phi_\textup{w}$ corresponds
to continuity (w.r.t.\ the domain) of the so-called optimal value function. 
Consulting \cite[Theorems~4.2.1, 4.2.2]{BankGuddatKlatteKummerTammer1983}, these properties are obtained
under continuity of the objective function and lower semicontinuity of the constraint mapping $\Gamma$.
In this regard, \cref{thm:closedness_of_graphs} amounts to a reasonable generalization to the
vector-valued situation. 

\begin{remark}\label{rem:nonclosedness_and_lack_of_usc_for_efficient_nondominated_map}
	Recall that the efficient set as well as the associated nondominated set in multiobjective optimization
	do not need to be closed even in situations where all objective functions are continuous and the feasible set is
	compact. This issue naturally extends to the parametric setting. 
	Typically, the set-valued mappings $\Psi$ and $\Phi$ do not possess closed graphs in many situations.
\end{remark}

Let us inspect some simple examples.

\begin{example} \label{ex:efficient_set_not_closed_paramindependent}
	We consider \eqref{eq:parametric_problem} with $n:=m:=1$, $q:=2$, $\mathcal C:=\R^2_+$, and
	\[
		\forall x,y\in\R\colon\quad
		f(x,y):=(\sin y,y),\qquad
		\Gamma(x) := [0,2 \pi].
	\]
	The efficiency map is given by 
	\begin{align*}
		\forall x \in \R \colon\quad 
		\Psi(x)= \{ 0 \} \cup (\pi, 3\pi/2]
	\end{align*}
	which is not closed for each $x \in \R$, whereas the weak efficiency map, given by
	\begin{align*}
		\forall x \in \R \colon\quad
		\Psi_{\textup{w}}(x)  = \{ 0 \} \cup [\pi, 3\pi/2] ,
	\end{align*}
	is closed at each $x \in \R$.  
\end{example}

\begin{example} \label{ex:efficiency_map_no_closed_graph}
	We consider \eqref{eq:parametric_problem} with $n:=1$, $m:=2$, $q:=2$, $\mathcal C:=\R^2_+$, and
	\[
		\forall x\in\R\,\forall y\in\R^2\colon\quad
		f(x,y):= y,
		\qquad
		\Gamma(x):=\{y \,|\,y_1 + x y_2 \geq 0\}.
	\]
	For $x > 0$, the efficient and weakly efficient points of \eqref{eq:parametric_problem} coincide. 
	In the case $x = 0$, there are no efficient point of \eqref{eq:parametric_problem}, 
	but all points from $\{0\}\times\R$ are weakly efficient. 
	Moreover, for $x < 0$, there is neither an efficient nor a weakly efficient point. 
	Hence, the efficiency map and frontier map are given by
	\begin{align*}
		\forall x \in \R \colon \quad
		\Psi(x) = \Phi(x) 
		= 
		\begin{cases}
			\{y \,|\, y_1 + x y_2 = 0\} & x > 0, \\
			\emptyset & x \leq 0,
		\end{cases}  
	\end{align*}
	and these mappings do not possess closed graphs. 
	Furthermore, the weak efficiency map and the weak frontier map are given by
	\begin{align*}
		\forall x \in \R \colon \quad
		\Psi_{\textup{w}}(x) = \Phi_{\textup{w}}(x) 
		= 
		\begin{cases}
			\{y \,|\, y_1 + x y_2 = 0\} & x \geq 0, \\
			\emptyset & x < 0
		\end{cases} 
	\end{align*}
	and possess closed graphs, see \Cref{cor:closedness_of_graphs} while noting that $\Gamma$ is
	not locally bounded at each point $x\geq 0$. 
\end{example}

\begin{remark}\label{rem:coincidence_of_efficiency_and_weak_efficiency}
	Assume that, for each $x\in\dom\Gamma$, $f(x,\cdot)\colon\R^m\to\R^q$ is $\mathcal C$-strictly-quasiconvex 
	while $\Gamma(x)$ is convex.
	Then, due to \cref{lem:closedness_of_efficient_set}, we find $\Psi(x)=\Psi_\textup{w}(x)$ and
	$\Phi(x)=\Phi_\textup{w}(x)$ for all $x\in\dom\Gamma$, i.e., there is no difference between
	the models \eqref{eq:BPP_eff} and \eqref{eq:BPP_weff}, 
	\eqref{eq:VFR_eff} and \eqref{eq:VFR_weff}, or
	\eqref{eq:VFR_eff_C} and \eqref{eq:VFR_weff_C}, respectively.
	Particularly, \cref{lem:closedness_of_efficient_set} and \cref{cor:closedness_of_graphs} can be
	used to infer closedness results on $\Psi$ and $\Phi$.
\end{remark}  

In the following lemma, we show that the efficiency and frontier mapping of a fully linear multiobjective
parametric optimization problem possess closed graphs.

\begin{lemma}\label{lem:closedness_in_linear_setting}
	For matrices $A\in\R^{r\times n}$, $B\in\R^{r\times m}$, $D\in\R^{q\times m}$, and $d\in\R^{r}$,
	let $f\colon\R^n\times\R^m\to\R^q$ and $\Gamma\colon\R^n\tto\R^m$ be given by
	\begin{equation}\label{eq:lower_level_linear}
		\forall x\in\R^n\,\forall y\in\R^m\colon\quad
		f(x,y):=Dy,\qquad
		\Gamma(x):=\{y\,|\,Ax+By\leq d\}.
	\end{equation}
	Furthermore, let the ordering cone $\mathcal C:=\R^q_+$ be fixed.
	Then the associated mappings $\Psi$ and $\Phi$ possess closed graphs.
\end{lemma}
\begin{proof}
	We prove the closedness of $\gph\Psi$ by combining a classical scalarization argument 
	and results from scalar linear parametric optimization.
	Afterwards, the closedness of $\gph\Phi$ is shown via Hofmann's lemma.
	
	Fix a sequence $((x^k,y^k))_{k\in\N}\subset\gph\Psi$ and a pair $(\bar x,\bar y)\in\R^n\times\R^m$
	such that $(x^k,y^k)\to(\bar x,\bar y)$.
	Let $k\in\N$ be fixed for the moment.
	Due to \cite[Theorem~4.7]{Ehrgott2005}, $y^k$ is an optimal
	solution of the scalar linear optimization problem
	\begin{equation}\label{eq:hybrid_scalarization}
		\min\limits_{y}\left\{\mathtt e^\top Dy\,\middle|\,Ax^k+By\leq d,\,Dy\leq Dy^k\right\}
	\end{equation}
	where $\mathtt e\in\R^q$ denotes the all-ones vector.
	The dual of this problem is given by
	\begin{equation}\label{eq:hybrid_scalarization_dual}
		\max\limits_{u,v}
			\left\{
				(d-Ax^k)^\top u+(Dy^k)^\top v
				\,\middle|\,
				B^\top u+D^\top v=D^\top\mathtt e,\,u\leq 0,\,v\leq 0
			\right\}.
	\end{equation}
	As \eqref{eq:hybrid_scalarization} possesses a solution, 
	\eqref{eq:hybrid_scalarization_dual} is solvable as well.
	
	Let us inspect \eqref{eq:hybrid_scalarization_dual} which can be interpreted as a
	linear parametric optimization problem which multiple parameters, namely $x^k$ and $y^k$,
	which merely appear in the objective function.
	Based on \cite[Sections~7, 8]{NozickaGuddatHollatzBank1974}, the set of parameters where
	this parametric problem possesses a solution is a convex polyhedron and can be partitioned into
	finitely many so-called regions of stability, being convex polyhedra as well, 
	such that, in the interior of any such region,
	precisely one of the vertices associated with the convex polyhedron
	\[
		Q:=\{(u,v)\,|\,B^\top u+D^\top v=D^\top\mathtt e,\,u\leq 0,\,v\leq 0\}
	\]
	is optimal, while on the boundary, this vertex is still
	optimal but may not be the unique maximizer.
	Anyhow, since there exist only finitely many regions of stability, there is an infinite index set $N\subset\N$
	such that all points from the subsequence $((x^k,y^k))_{k\in N}$ belong
	to the same region of stability which is associated 
	with some fixed vertex $(\bar u,\bar v)\in\R^n\times\R^m$ of $Q$.
	As all regions of stability are closed, $(\bar u,\bar v)$ does not only solve
	\eqref{eq:hybrid_scalarization_dual} for each $k\in N$ but also the problem 
	\begin{equation}\label{eq:hybrid_scalarization_dual_limit}
		\max\limits_{u,v}
			\left\{
				(d-A\bar x)^\top u+(D\bar y)^\top v
				\,\middle|\,
				(u,v)\in Q
			\right\}.
	\end{equation}
	Strong duality of linear programming yields $\mathtt e^\top Dy^k=(d-Ax^k)^\top \bar u+(Dy^k)^\top\bar v$
	for each $k\in N$, and taking the limit $k\to_N\infty$ gives 
	$\mathtt e^\top D\bar y=(d-A\bar x)^\top\bar u+(D\bar y)^\top\bar v$.
	As $\bar y$ is feasible to 
	\[
		\min\limits_{y}\left\{\mathtt e^\top Dy\,\middle|\,A\bar x+By\leq d,\,Dy\leq D\bar y\right\},
	\]
	which is the dual of \eqref{eq:hybrid_scalarization_dual_limit},
	it already must be a minimizer.
	Applying \cite[Theorem~4.7]{Ehrgott2005} once more gives $\bar y\in\Psi(\bar x)$.
	This shows the closedness of $\gph\Psi$.
		
	In the remainder of the proof, we show that $\gph\Phi$ is closed.
	Thus, let $(\bar x,\bar z)\in\R^n\times\R^q$ and $((x^k,z^k))_{k\in\N}\subset\gph\Phi$ be chosen
	such that $(x^k,z^k)\to(\bar x,\bar z)$.
	By definition of $\Phi$, we find a sequence $(y^k)_{k\in\N}\subset\R^m$ such that
	$y^k\in\Psi(x^k)$ and $z^k=Dy^k$ hold for each $k\in\N$.
	Applying Hofmann's lemma, see e.g.\ \cite[Lemma~3C.4]{DontchevRockafellar2014}, for each $k\in\N$
	to the linear system
	\[
		Dy\leq z^k,\,-Dy\leq -z^k,\,By\leq d-Ax^k
	\]
	yields some $\tilde y^k\in\Gamma(x^k)$ which satisfies $D\tilde y^k=z^k$ and
	$\nnorm{\tilde y^k}\leq \kappa(\nnorm{z^k}+\nnorm{d-Ax^k})$
	where the constant $\kappa>0$ does not depend on $k$.
	Particularly, for each $k\in\N$, $D\tilde y^k=z^k=Dy^k$, i.e., $\tilde y^k\in\Psi(x^k)$.
	Further, by boundedness of $(x^k)_{k\in\N}$ and $(z^k)_{k\in\N}$, $(\tilde y^k)_{k\in\N}$ is bounded as well so we
	may assume $\tilde y^k\to\tilde y$ for some $\tilde y\in\R^m$.
	The first part of the proof gives $\tilde y\in\Psi(\bar x)$, so that $D\tilde y=\bar z$ yields $\bar z\in\Phi(\bar x)$.
	Hence, $\gph\Phi$ is closed.
\end{proof}

We note that the first part of the above proof has been inspired by classical arguments of Isermann
who showed that all efficient points of a linear multiobjective optimization problem can be found
with the aid of the weighted-sum scalarization technique where all weights are positive,
see \cite{Isermann1974}.

In the context of \cref{prop:value_function_char_extended_of_eff_set} and \cref{cor:adding_an_ordering_cone_in_VFR},
the following result should be taken into account.
\begin{lemma}\label{lem:closedness_in_C_extended_approach}
	Let $\Gamma$ be locally bounded at each point of its domain.
	Then the following assertions hold.
	\begin{enumerate}
		\item The set $\gph\Phi$ is closed if and only if $\gph(\Phi-\mathcal C)$ is closed.
		\item The set $\gph\Phi_\textup{w}$ is closed if and only if $\gph(\Phi_\textup{w}-\mathcal C)$ is closed.
	\end{enumerate}
\end{lemma}
\begin{proof}
	We only prove the first assertion. 
	The second one can be shown in analogous fashion 
	while incorporating some arguments from the proof of \cref{prop:value_function_char_extended_of_eff_set}.
	
	We show both implications separately.\\
	$[\Longrightarrow]$ Let $\gph\Phi$ be closed, 
	and choose a sequence $((x^k,z^k))_{k\in\N}\subset\gph(\Phi-\mathcal C)$ such that
	$(x^k,z^k)\to(\bar x,\bar z)$ for some $(\bar x,\bar z)\in\R^n\times\R^q$.
	By definition, we find sequences $(c^k)_{k\in\N}\subset\mathcal C$ and $(y^k)_{k\in\N}\subset\R^m$ such that
	$z^k=f(x^k,y^k)-c^k$ and $y^k\in\Psi(x^k)$, i.e., $(x^k,f(x^k,y^k))\in\gph\Phi$, for each $k\in\N$.
	As the sequence $((x^k,y^k))_{k\in\N}\subset\gph\Gamma$ is bounded by assumption, it possesses an accumulation
	point $(\bar x,\bar y)$ which, by continuity of $f$ and closedness of $\gph\Phi$, 
	satisfies $f(\bar x,\bar y)\in\Phi(\bar x)$.
	Furthermore, $(c^k)_{k\in\N}$ possesses the accumulation point $f(\bar x,\bar y)-\bar z$ which, by closedness
	of $\mathcal C$, belongs to $\mathcal C$. Hence, we have shown $\bar z\in(\Phi-\mathcal C)(\bar x)$, and
	closedness of $\gph(\Phi-\mathcal C)$ follows.\\
	$[\Longleftarrow]$ Let $\gph(\Phi-\mathcal C)$ be closed,
	and choose a sequence $((x^k,z^k))_{k\in\N}\subset\gph\Phi$ such that
	$(x^k,z^k)\to(\bar x,\bar z)$ for some $(\bar x,\bar z)\in\R^n\times\R^q$.
	On the one hand, since $0\in\mathcal C$, we also have $((x^k,z^k))_{k\in\N}\subset\gph(\Phi-\mathcal C)$, and
	$(\bar x,\bar z)\in\gph(\Phi-\mathcal C)$ follows, i.e., there are $\bar y\in\Psi(\bar x)$ and $c\in\mathcal C$
	such that $\bar z=f(\bar x,\bar y)-c$.
	On the other hand, there is a sequence $(y^k)_{k\in\N}\subset\R^m$ such that $z^k=f(x^k,y^k)$ and
	$y^k\in\Psi(x^k)$ for all $k\in\N$. 
	As the sequence $((x^k,y^k))_{k\in\N}\subset\gph\Gamma$ is bounded by assumption, it possesses an
	accumulation point $(\bar x,\tilde y)$ which, by closedness of $\gph\Gamma$, belongs to $\gph\Gamma$.
	Combining these observations, we have $f(\bar x,\tilde y)=\bar z=f(\bar x,\bar y)-c$.
	Now, $\bar y\in\Psi(\bar x)$ guarantees $c=0$, i.e., $\bar z=f(\bar x,\bar y)\in\Phi(\bar x)$,
	and closedness of $\gph\Phi$ follows.
\end{proof}

\subsection{Existence results in multiobjective bilevel optimization}\label{sec:existence}

The following example shows that the closedness of $\gph \widehat \Psi$ is, indeed,
essential for the existence of efficient and weakly efficient points of \eqref{eq:BPP}. 

\begin{example} \label{ex:efficient_set_not_closed}
	We consider \eqref{eq:BPP} with $n := m := 1$, $p := 2$, $\mathcal K := \R^2_+$, $X:=\R$, and
	\[
		\forall x,y\in\R\colon\quad
		F(x,y):=(\cos(y) + 1,(y - \pi)^2).
	\]
	Furthermore, we assume that $\hat\Psi \in \{ \Psi, \Psi_{\textup{w}} \} $ is the efficiency or weak efficiency map of the parametric optimization
	problem from \cref{ex:efficient_set_not_closed_paramindependent}, i.e.,
	\[
		\forall x \in \R \colon \quad
		\Psi(x) = \{ 0 \} \cup (\pi, 3\pi/2],
		\qquad
		\Psi_\textup{w}(x) = \{0\} \cup [\pi,3\pi/2]. 
	\]
	In the case $\widehat\Psi := \Psi$, there is neither an efficient nor a weakly effficient point of \eqref{eq:BPP}.
	However, in the case  $\widehat\Psi := \Psi_{\textup{w}}$, $\R\times\{\pi\}$ is a feasible subset of \eqref{eq:BPP},
	and these are precisely the efficient and, simultaneously, weakly efficient points of the problem. 
\end{example}

The following result concerns
the existence of efficient and weakly efficient points of \eqref{eq:BPP}. 

\begin{theorem} \label{thm:existence_of_efficient_points_of_BPP}
	Let $\gph\widehat\Psi$ be closed and $\gph\widehat\Psi \cap (X \times \R^m)$ be nonempty and bounded. 
	Then there exist an efficient and a weakly efficient point of \eqref{eq:BPP}.
\end{theorem}
\begin{proof}
	Problem \eqref{eq:BPP} is equivalent to 
	\begin{equation}\label{eq:BPP_graph_Psi}
		\min\limits_{x,y}{}_{\mathcal{K}}\{F(x,y)\,|\,(x,y) \in \gph\widehat\Psi\cap(X \times \R^m) \}.
	\end{equation}
	Due to the postulated assumptions, 
	the set $\gph\widehat\Psi \cap (X \times \R^m)$ is nonempty, bounded, and closed and, thus, compact. 
	Then $F(\gph\widehat\Psi \cap (X \times \R^m))$ is compact because $F$ is continuous. 
	The existence of an efficient point of \eqref{eq:BPP_graph_Psi}, 
	which then is efficient and, thus, weakly efficient for \eqref{eq:BPP}, 
	follows by \cite[Theorem~6.3]{Jahn2011}.
\end{proof}

Let us note that the boundedness assumption on $\gph\widehat\Psi\cap(X\times\R^m)$ in \cref{thm:existence_of_efficient_points_of_BPP} 
is inherently satisfied whenever $X$ and $\bigcup_{x\in X}\Gamma(x)$ are bounded.

Based on \cref{cor:closedness_of_graphs}, 
we obtain from \Cref{thm:existence_of_efficient_points_of_BPP} an existence result for efficient and weakly efficient points of \eqref{eq:BPP_weff}.

\begin{corollary}\label{cor:existence_weakly_efficient_points_lower_level}
	Let $\dom\Gamma$ be closed, and let $\Gamma$ be locally bounded and lower semicontinuous
	w.r.t.\ its domain at each point from $\dom\Gamma$. 
	Let $\gph\Psi_\textup{w} \cap (X \times \R^m)$ be nonempty and bounded. 
	Then there exist an efficient and a weakly efficient point of \eqref{eq:BPP_weff}.
\end{corollary}

We note that \cref{cor:existence_weakly_efficient_points_lower_level} generalizes \cite[Theorem~5.5]{DempeMehlitz2019} 
to multiobjective bilevel optimization problems, but does not cover the situation were the lower-level variable comes
from an infinite-dimensional space.

Similar as in \cref{cor:existence_weakly_efficient_points_lower_level}, 
\cref{thm:existence_of_efficient_points_of_BPP} gives an existence result for efficient and weakly efficient points
of \eqref{eq:BPP_eff} and  \eqref{eq:BPP_weff} 
whenever the lower-level data is fully linear as described in \eqref{eq:lower_level_linear}.

\begin{corollary}
	Let the lower-level data be given as in \eqref{eq:lower_level_linear} and fix $\mathcal C:=\R^q_+$.
	Then the following assertions hold.
	\begin{enumerate}
		\item Let $\gph\Psi \cap (X \times \R^m)$ be nonempty and bounded. 
			Then there exist an efficient and a weakly efficient point of \eqref{eq:BPP_eff}.
		\item Let $\gph\Psi_\textup{w}\cap(X\times\R^m)$ be nonempty and bounded.
			Then there exist an efficient and a weakly efficient point of \eqref{eq:BPP_weff}.
	\end{enumerate}
\end{corollary}
\begin{proof}
	The first assertion directly follows from \cref{lem:closedness_in_linear_setting} and \cref{thm:existence_of_efficient_points_of_BPP}.
	The second statement can be obtained from \cref{cor:closedness_of_graphs} and \cref{thm:existence_of_efficient_points_of_BPP}
	as $\gph\Gamma$ is a convex polyhedral set and, thus, closed, and these properties also ensure that $\Gamma$ is lower
	semicontinuous w.r.t.\ its domain at each point from $\dom\Gamma$, see e.g.\ \cite[Example~9.35]{RockafellarWets2009}.
\end{proof}

\begin{remark}\label{rem:weaker_closedness_assumptions_for_existence}
	Taking a look back at the optimal value function reformulation for standard scalar
	bilevel optimization problems, 
	we know from \cref{lem:closed_graph_for_special_map}
	that the feasible set of \eqref{eq:VFR_classic_with_ordering_cone} 
	is closed if and only if $\varphi$ is
	upper semicontinuous, and this can also be seen from \eqref{eq:VFR_classic}
	as sublevel sets of extended-real-valued functions are closed if and only if
	they are lower semicontinuous.
	The latter property is clearly milder than the continuity of $\varphi$ 
	on its domain $\{x\,|\,|\varphi(x)|<\infty\}$ which seems to
	be necessary for this desirable closedness when focusing on the model problem
	\eqref{eq:VFR_nonstandard}.
	
	Thus, having problems \eqref{eq:VFR_eff_C} and \eqref{eq:VFR_weff_C} 
	as well as \cref{cor:adding_an_ordering_cone_in_VFR} at hand, 
	it is seemingly reasonable to work out conditions which ensure that the graphs of
	$\Phi-\mathcal C$ and $\Phi_\textup{w}-\mathcal C$ are closed.
	This is enough to obtain
	the closedness of the feasible sets of \eqref{eq:VFR_eff_C} and \eqref{eq:VFR_weff_C}
	which, in turn, provides the basis for further existence results in
	multiobjective bilevel optimization.
	However, note that due to \cref{lem:closedness_in_C_extended_approach}, 
	it turns out that whenever the feasibility mapping $\Gamma$ is locally bounded,
	then this apparently milder approach yields no new insights.
	Observe that this is in accordance with the aforementioned scalar situation as
	the optimal value function $\varphi$ is lower semicontinuous on its domain 
	whenever $\Gamma$ is locally bounded.
\end{remark}

\subsection{An intermediate approach to multiobjective bilevel optimization}\label{sec:intermediate_approach}

On the one hand, the feasible set of \eqref{eq:BPP_eff} is likely to be much smaller than the one of
\eqref{eq:BPP_weff} since the solution criterion behind $\Psi$ is sharper than 
the one behind $\Psi_\textup{w}$.
On the other hand, from the viewpoint of existence theory, necessary optimality conditions, and solution
algorithms, a closed feasible set is desirable. The latter is given for \eqref{eq:BPP_weff} under mild
assumptions, see \cref{cor:closedness_of_graphs}, while the feasible set of \eqref{eq:BPP_eff} is not
likely to be closed apart from selected situations.

In this subsection, we suggest yet another possibility for the choice of $\widehat\Psi$ 
in the general multiobjective bilevel optimization problem \eqref{eq:BPP} 
which induces an intermediate model between \eqref{eq:BPP_eff} and \eqref{eq:BPP_weff}. 
To start, however, we would like to suggest an intermediate value function model
which bridges \eqref{eq:VFR_eff} and \eqref{eq:VFR_weff} to some extent.
There are two reasons for this approach. 
First, it allows for a simple definition of an intermediate solution map such that the
associated models \eqref{eq:BPP} and \eqref{eq:value_function_reformulation} are
equivalent. Second, noting that, in order to determine efficient and weakly efficient
points of some multiobjective optimization problem, one has to compute its nondominated and
weakly nondominated points first before identifying their pre-images in the feasible set,
this procedure seems the more natural one.

Therefore, let us introduce $\overline\Phi\colon\R^n\tto\R^q$ by means of
\[
	\gph\overline\Phi:=\cl\gph\Phi
\]
and consider
\begin{equation}\label{eq:VFR_inter_eff}\tag{$\overline{\mathcal E}$-VFR}
	\min\limits_{x,y}{}_{\mathcal{K}}\{F(x,y)\,|\,x\in X,\,f(x,y)\in\overline\Phi(x),\,y\in\Gamma(x)\}.
\end{equation}
As $\gph\Phi_\textup{w}$ is closed under mild assumptions, see e.g.\ \cref{cor:closedness_of_graphs},
the inclusions $\gph\Phi\subset\gph\overline\Phi\subset\gph\Phi_\textup{w}$ are likely to hold
and would give $\Phi(x)\subset\overline\Phi(x)\subset\Phi_\textup{w}(x)$ for each $x\in\R^n$.
Hence, \eqref{eq:VFR_inter_eff} can be interpreted as an intermediate problem between
\eqref{eq:VFR_eff} and \eqref{eq:VFR_weff} since the feasible set of \eqref{eq:VFR_inter_eff}
is not smaller than the one of \eqref{eq:VFR_eff} and likely to be not larger than the one of
\eqref{eq:VFR_weff}. 
Furthermore, \eqref{eq:VFR_inter_eff} possesses a closed feasible set by construction of $\overline\Phi$,
and, thus, efficient and weakly efficient points of \eqref{eq:VFR_inter_eff} exist under standard assumptions.

Based on $\overline\Phi$, let us introduce $\overline\Psi\colon\R^n\tto\R^m$ by means of
\[
	\forall x\in\R^n\colon\quad
	\overline\Psi(x):=\{y\in\Gamma(x)\,|\,f(x,y)\in\overline\Phi(x)\}.
\]
In the subsequent lemma, we comment on the properties of $\overline\Psi$.
The simple proof follows by definition of $\overline\Psi$ as well as our above discussion on the
properties of $\overline\Phi$ and is, thus, omitted.
\begin{lemma}\label{lem:solution_map_intermediate_approach}
	\leavevmode
	\begin{enumerate}
		\item The set $\gph\overline\Psi$ is closed.
		\item\label{item:inclusions_for_Psi_maps} 
			Whenever $\gph\Phi_\textup{w}$ is closed, 
			then the inclusions $\gph\Psi\subset\gph\overline\Psi\subset\gph\Psi_\textup{w}$ hold.
			Particularly, $\Psi(x)\subset\overline\Psi(x)\subset\Psi_\textup{w}(x)$ is valid for
			each $x\in\R^n$ in this situation.
		\item\label{item:closure_of_graph_only_subset} 
			We have $\cl\gph\Psi\subset\gph\overline\Psi$.
	\end{enumerate}
\end{lemma}

By definition of $\overline\Psi$, the problem
\begin{equation}\label{eq:BPP_int_eff}\tag{$\overline{\mathcal{E}}$-BPP}
	\min\limits_{x,y}{}_{\mathcal{K}}\{F(x,y)\,|\,x\in X,\,y\in\overline{\Psi}(x)\}
\end{equation}
is equivalent to \eqref{eq:VFR_inter_eff}, 
and due to \cref{lem:solution_map_intermediate_approach}~\ref{item:inclusions_for_Psi_maps},
can be seen as an intermediate model between \eqref{eq:BPP_eff} and \eqref{eq:BPP_weff}.
Clearly, \eqref{eq:BPP_int_eff} has a closed feasible set and, 
due to \cref{thm:existence_of_efficient_points_of_BPP},
is likely to possess efficient and weakly efficient points which can be interpreted as
approximately efficient and approximately weakly efficient for \eqref{eq:BPP_eff}.

	Let us state yet another motivation behind the consideration of the models
	\eqref{eq:BPP_int_eff} and \eqref{eq:VFR_inter_eff}.
	When solving \eqref{eq:BPP_eff} and \eqref{eq:VFR_eff} numerically with a suitable
	solution algorithm, respectively,
	one typically obtains a sequence of iterates from the solution method.
	Taking into account the potential non-closedness of $\gph\Psi$ and $\gph\Phi$,
	accumulation points of this sequence are not necessarily feasible to
	\eqref{eq:BPP_eff} and \eqref{eq:VFR_eff}, respectively.
	However, they are likely to be feasible to
	\eqref{eq:BPP_int_eff} and \eqref{eq:VFR_inter_eff}
	as $\gph\overline\Psi$ and $\gph\overline\Phi$ are closed
	supersets of $\gph\Psi$ and $\gph\Phi$, respectively.
	Hence, qualitative properties of solution methods associated with
	\eqref{eq:BPP_eff} and \eqref{eq:VFR_eff} are naturally related to the model problems
	\eqref{eq:BPP_int_eff} and \eqref{eq:VFR_inter_eff}, respectively.

Let us point out that one does not obtain equality in 
\cref{lem:solution_map_intermediate_approach}~\ref{item:closure_of_graph_only_subset}
in general situations. Furthermore, when defining $\widetilde\Psi\colon\R^n\tto\R^m$ by
$\gph\widetilde\Psi:=\cl\gph\Psi$ and $\widetilde\Phi\colon\R^n\tto\R^q$ by
$\widetilde\Phi(x):=f(x,\widetilde\Psi(x))$ for all $x\in\R^n$, 
then the corresponding problems \eqref{eq:BPP} and \eqref{eq:value_function_reformulation}
for $\widehat\Psi:=\widetilde\Psi$ and $\widehat\Phi:=\widetilde\Phi$, respectively, are not
necessarily equivalent. We omitted this undesirable behavior by starting with the
definition of a suitable intermediate frontier mapping before introducing an
intermediate efficiency mapping.

\begin{example} \label{ex:order_of_definitions_for_intermediate_approach_is_important}
We consider \eqref{eq:parametric_problem} with $n:=1$, $m:=q:=2$, $\mathcal C:= \R^2_+$, and
\begin{align*}
	\forall x \in \R\,\forall y\in\R^2\colon\quad
	f(x,y) 		:= (\sin y_1, y_2),\qquad
	\Gamma(x) 	:= \conv\{(0,0),(2\pi,2\pi)\} \cup \{ (0,\pi) \} 
\end{align*}
We easily find
\begin{align*} 
	\forall x\in \R\colon\quad
	\Psi(x) &= \{ y \, | \, y_1 = y_2, y_2 \in \{ 0 \} \cup (\pi, 3\pi/2] \},\\
	\Phi(x) &= \{(\sin z,z)\,|\,z\in\{0\}\cup(\pi,3\pi/2]\}
\end{align*}
which gives 
\begin{align*}
	\forall x \in \R\colon\quad
	\widetilde\Psi(x) &= \{ y \, | \, y_1 = y_2, y_2 \in \{ 0 \} \cup [\pi, 3\pi/2] \},\\
	\widetilde\Phi(x) &= \{ (\sin z, z) \, | \, z \in \{ 0 \} \cup [\pi, 3\pi/2] \}. 
\end{align*}
Observe, however, that
\[ \forall x \in \R\colon\quad
\{ y \in \Gamma(x) \, | \; f(x,y) \in \widetilde\Phi(x) \} = \widetilde\Psi(x) \cup \{ (0,\pi) \},
\]
i.e., the feasible set of any superordinate model \eqref{eq:BPP} for $\widehat\Psi:=\widetilde\Psi$ in likely
to be smaller than the feasible set of the associated value function reformulation 
\eqref{eq:value_function_reformulation} for $\widehat\Phi:=\widetilde\Phi$.

Let us also mention that
\begin{align*} 
	\forall x \in \R\colon\quad
		\overline\Phi(x) = \{ (\sin z, z) \, | \, z \in \{ 0 \} \cup [\pi, 3\pi/2] \} = \widetilde\Phi(x),\qquad
		\overline{\Psi}(x) = \widetilde\Psi(x) \cup \{ (0,\pi) \}
\end{align*}
holds in this example. This also shows $\cl\gph\Psi\subsetneq\gph\overline\Psi$.
\end{example}

The following example shows that 
\eqref{eq:BPP_eff}, \eqref{eq:BPP_weff}, and \eqref{eq:BPP_int_eff}
may possess pairwise different solutions.
\begin{example}\label{ex:intermediate_approach}
	Motivated by \cref{ex:domination_property_without_closedness}, for $n:=1$, $m:=q:=2$, and $\mathcal C:=\R^2_+$,
	we consider the parametric multiobjective optimization problem 
	\eqref{eq:parametric_problem} with the data
	\[
		\forall x\in\R\,\forall y\in\R^2\colon\quad
		f(x,y):=y,\quad
		\Gamma(x):=\{y\in\R^2_+\,|\,y_1^2+y_2^2\geq x^2\}\cup([-1,0]\times[|x|,\infty)).
	\]
	A simple calculation reveals
	\begin{align*}
		\forall x\in\R\colon\quad
		\Psi(x)&=
			\begin{cases}
				\{(|x|\,\cos t,|x|\,\sin t)\,|\,t\in[0,\pi/2)\}\cup\{(-1,|x|)\}	&	x\neq 0,\\
				\{(-1,0)\}	&	x=0,
			\end{cases}
			\\
		\Psi_\textup{w}(x)	&=
			\{(|x|\,\cos t,|x|\,\sin t)\,|\,t\in[0,\pi/2]\}\cup([-1,0]\times\{|x|\})\\
			&\qquad
			\cup([|x|,\infty)\times\{0\})\cup(\{-1\}\times[|x|,\infty)),\\
		\overline{\Psi}(x) &=
			\{(|x|\,\cos t,|x|\,\sin t)\,|\,t\in[0,\pi/2]\}\cup\{(-1,|x|)\}.
	\end{align*}
	Observe that $\gph\overline{\Psi}=\gph\Psi\cup\{(x,(0,|x|))\,|\,x\in\R\}$ holds true.
	
	Now, we consider the superordinate semivectorial bilevel optimization problem
	\[
		\min\limits_{x,y}\{(y_1+1/4)^2+y_2^2\,|\,x\geq 1,\,y\in\widehat\Psi(x)\}
	\]
	where $\widehat\Psi\in\{\Psi,\Psi_\textup{w},\overline\Psi\}$.
	Let us start with the investigation of $\widehat\Psi:=\Psi$.
	The associated problem \eqref{eq:BPP_eff} does not possess a global minimizer,
	but there is a local minimizer at $(1,(-1,1))$.
	Whenever $\widehat\Psi:=\Psi_\textup{w}$, the associated problem \eqref{eq:BPP_weff}
	possesses the uniquely determined global minimizer $(1,(-1/4,1))$, and no other
	local minimizers exist.
	Finally, for $\widehat\Psi:=\overline{\Psi}$, the associated problem \eqref{eq:BPP_int_eff}
	possesses the uniquely determines global minimizer $(1,(0,1))$ 
	(which is not feasible to \eqref{eq:BPP_eff}) 
	and an additional local minimizer at $(1,(-1,1))$.
\end{example}

Finally, we would like to point the reader's attention to the fact that, 
in contrast to \cref{cor:adding_an_ordering_cone_in_VFR},
\eqref{eq:VFR_inter_eff} and
\begin{equation}\label{eq:VFR_inter_eff_C}
	\min\limits_{x,y}{}_{\mathcal{K}}\{F(x,y)\,|\,x\in X,\,f(x,y)\in\overline\Phi(x)-\mathcal C,\,y\in\Gamma(x)\}
\end{equation}
are not necessarily equivalent as the feasible set of
\eqref{eq:VFR_inter_eff_C} is likely to be larger than the one
of \eqref{eq:VFR_inter_eff}.
However, if $\gph\Phi_\textup{w}$ is closed, the feasible set of
\eqref{eq:VFR_inter_eff_C} is included in the feasible set of
\eqref{eq:VFR_weff_C} which, by \cref{cor:adding_an_ordering_cone_in_VFR},
equals the feasible set of \eqref{eq:VFR_weff}. 
Thus, the enlargement of the feasible set in \eqref{eq:VFR_inter_eff_C}
in comparison with \eqref{eq:VFR_inter_eff} is somewhat controllable.

\begin{example}\label{eq:VFR_inter_eff_C_not_reasonable}
	Let us reinspect the semivectorial bilevel optimization problem from
	\cref{ex:intermediate_approach}.
	One can easily check that, for $\bar x:=1$, $\bar y:=(-1/4,1)\in\Gamma(\bar x)$ satisfies
	\[
		f(\bar x,\bar y)\in\overline\Phi(\bar x)-\R^2_+=\overline\Psi(\bar x)-\R^2_+,
	\]
	but $(\bar x,\bar y)$ is not feasible to the associated problem
	\eqref{eq:BPP_int_eff} and, thus, \eqref{eq:VFR_inter_eff}.
	Keeping our above discussion in mind and noticing that $\gph\Phi_\textup{w}$ is closed,
	the associated model \eqref{eq:VFR_inter_eff_C} possesses the global minimizer
	$(\bar x,\bar y)$, see \cref{ex:intermediate_approach} as well.
\end{example}

\section{Towards necessary optimality conditions}\label{sec:coderivative_calculus}

For the derivation of necessary optimality conditions and constraint qualifications for 
\eqref{eq:BPP_eff} and \eqref{eq:BPP_weff} via \eqref{eq:VFR_eff} and \eqref{eq:VFR_weff},
respectively, in terms of the limiting tools of generalized differentiation, it is 
important to have formulas or at least upper estimates for the coderivative of
$\Phi$ and $\Phi_\textup{w}$ available.
Here, we focus on the latter and comment on two approaches which can be used to obtain
such estimates in terms of initial problem data.

\subsection{Coderivative estimates via the domination property}\label{sec:coderivative_domination}

In the literature, one can find several papers dealing with the derivation of upper
coderivative estimates for the frontier mapping $\Phi$ in the case where \eqref{eq:parametric_problem}
possesses the domination property for all $x\in\R^n$ from a neighborhood of the reference
point, see e.g.\ \cite{HuyMordukhovichYao2008,LiXue2014,XueLiLiaoYao2011}.
Apart from the fact that these papers completely ignore the fact that $\gph\Phi$ is likely to
be non-closed, so that the coderivative is not well-defined, they obtain validity of the
upper estimate
\begin{equation}\label{eq:upper_estimate_frontier_map}
	D^*\Phi(\bar x,\bar z)(z^*)\subset D^*(\Sigma+\mathcal C)(\bar x,\bar z)(z^*)
\end{equation}
only for points $(\bar x,\bar z)\in\gph\Phi$ and $z^*\in \mathcal{C}^*_>$ where
\[
	\mathcal{C}^*_>:=\{z^*\,|\,\forall z\in \mathcal C\cap\mathbb S_1(0)\colon\,\langle z^*,z\rangle>0\}.
\]
In practice, this means that the estimate \eqref{eq:upper_estimate_frontier_map} is, more or less,
only reasonable for arguments $z^*\in\R^q$ from $\intr \mathcal{C}^*$.
For the derivation of necessary optimality conditions for \eqref{eq:BPP_eff} via \eqref{eq:VFR_eff}, 
however, this additional constraint
on the variable $z^*$ might be too restrictive.

In the proof of \cite[Theorem~4.1]{LafhimZemkoho2022}, the authors claim that
for $\mathcal C:=\R^q_+$, one has the estimate
\begin{equation}\label{eq:wrong_estimate}
	D^*\Phi(\bar x,\bar z)(z^*)\subset D^*\Sigma(\bar x,\bar z)(z^*)
\end{equation}
for all $z^*\in\R^q$. This estimate, however, does not depend on the choice of the
ordering cone $\mathcal C$ which is somewhat dubious. 
As the following example shows, this estimate is indeed false, and so are
\cite[Theorems~4.1, 4.4]{LafhimZemkoho2022}.
\begin{example}\label{ex:coderivative_of_frontier_map}
	We consider \eqref{eq:parametric_problem} with $n:=1$, $m:=q:=2$, $\mathcal C:=\R^2_+$, and
	\[
		\forall x\in\R\,\forall y\in\R^2\colon\quad
		f(x,y):=y,\qquad
		\Gamma(x):=\{y\,|\,y_1+2y_2\geq 0,\,2y_1+y_2\geq 0,\,y_1,y_2\leq 2\},
	\]
	i.e., the associated problem \eqref{eq:parametric_problem} does not depend on the parameter.
	We find
	\[
		\forall x\in\R\colon\quad
		\Sigma(x)=\Gamma(x),\quad
		\Phi(x)=\Psi(x)=\conv\{(-1,2),(0,0)\}\cup \conv\{(2,-1),(0,0)\}
	\]
	and observe that \eqref{eq:parametric_problem} possesses the domination property (for each $x\in\R$).
	Let us fix $\bar x:=0$ and $\bar z:=(0,0)$.
	Then we find
	\begin{align*}
		N_{\gph\Sigma}(\bar x,\bar z)&=\cone\{(0,-1,-2),(0,-2,-1)\},\\
		N_{\gph\Phi}(\bar x,\bar z)&=\cone\{(0,-1,-2),(0,-2,-1)\}\cup\cone\{(0,2,1)\}\cup\cone\{(0,1,2)\}.
	\end{align*} 
	This can be easily seen from $\gph\Sigma=\R\times\Sigma(0)$ and $\gph\Phi=\R\times\Phi(0)$,
	the product rule for the computation of limiting normals, see e.g.\ \cite[Proposition~1.2]{Mordukhovich2006},
	and a simple calculation of the normal cones to $\Sigma(0)$ (which is convex) and $\Phi(0)$ in $\R^2$.
	Hence, for $z^*:=(-1,-2)$, we find
	\[
		D^*\Phi(\bar x,\bar z)(z^*)=\{0\},
		\qquad
		D^*\Sigma(\bar x,\bar z)(z^*)=\emptyset,
	\]
	showing that the estimate \eqref{eq:wrong_estimate} fails to be true for general arguments.
\end{example}

Our next example indicates that the results from 
\cite{HuyMordukhovichYao2008,LiXue2014,XueLiLiaoYao2011} can be directly generalized
to weak frontier mappings.
\begin{example}\label{ex:coderivative_of_weak_frontier_map}
	We consider \eqref{eq:parametric_problem} with $n:=m:=1$, $q:=2$, $\mathcal C:=\R^2_+$, and
	\[
		\forall x,y\in\R\colon\quad
		f(x,y):=(xy,y),\qquad
		\Gamma(x):=[0,1].
	\]
	Then we find $\Sigma(x)=\conv\{(0,0),(x,1)\}$ for all $x\in\R$ as well as
	\[
		\forall x\in\R\colon\quad
		\Phi_\textup{w}(x)
		=
		\begin{cases}
			\Sigma(x)	&	x\leq 0,\\
			\{(0,0)\}	&	x>0,
		\end{cases}
		\quad
		\Psi_\textup{w}(x)
		=
		\begin{cases}
			[0,1]	&	x\leq 0,\\
			\{0\}	&	x>0.
		\end{cases}
	\]
	Using this, one can easily check that $\Sigma(x)\subset\Phi_\textup{w}(x)+\R^2_+$ holds for
	all $x\in\R$, i.e., the associated problem \eqref{eq:parametric_problem} possesses the
	weak domination property for all $x\in\R$ (actually, it also possesses the domination
	property for all $x\in\R$).
	Let us consider $\bar x:=\bar y:=0$ and $\bar z:=(0,0)$.
	It is not hard to see that
	\begin{align*}
		N_{\gph\Sigma}(\bar x,\bar z)			&=\{0\}\times\R\times\R_-,\\
		N_{\gph(\Sigma+\mathcal C)}(\bar x,\bar z)		&=\{0\}\times\R_-\times\R_-,\\
		N_{\gph\Phi_\textup{w}}(\bar x,\bar z)	&=(\{0\}\times\R\times\R)\cup(\R_+\times\R\times\{0\}),\\
		N_{\gph\Psi_\textup{w}}(\bar x,\bar y)	&=(\{0\}\times\R)\cup(\R_+\times\{0\}).
	\end{align*}
	This directly shows that the inclusion
	\begin{equation}\label{eq:upper_estimate_weak_frontier_map}
		D^*\Phi_\textup{w}(\bar x,\bar z)(z^*)\subset D^*(\Sigma+\mathcal C)(\bar x,\bar z)(z^*)
	\end{equation}
	holds for all $z^*\in\intr\R^2_+$, but is violated for, e.g., $z^*:=(1,0)$ already, where
	\[
		D^*\Phi_\textup{w}(\bar x,\bar z)(z^*)=\R_+,
		\qquad
		D^*(\Sigma+\mathcal C)(\bar x,\bar z)(z^*)=\{0\}.
	\]
\end{example}

Let us now formally prove that the estimate \eqref{eq:upper_estimate_weak_frontier_map} holds
for all $z^*\in \mathcal{C}^*_>$.

\begin{lemma}\label{lem:coderivative_estimates}
		Let $(\bar x,\bar z)\in\gph\Phi_\textup{w}$ be fixed and let $\gph\Phi_\textup{w}$ and
		$\gph(\Sigma+\mathcal C)$ be locally closed around $(\bar x,\bar z)$.
		Assume that \eqref{eq:parametric_problem} possesses the weak domination property for all $x\in U$
		where $U\subset\R^n$ is a neighborhood of $\bar x$.
		Then, for each $z^*\in \mathcal{C}^*_>$, the estimate \eqref{eq:upper_estimate_weak_frontier_map}
		is valid.
\end{lemma}
\begin{proof}
		Fix $z^*\in \mathcal{C}^*_>$.
		By mimicking the proof of \cite[Lemma~3.2]{HuyMordukhovichYao2008}, which does not
		rely on the precise underlying notion of nondominance, we find
		\[
			D^*\Phi_\textup{w}(\bar x,\bar z)(z^*)
			\subset
			D^*(\Phi_\textup{w}+\mathcal C)(\bar x,\bar z)(z^*).
		\]
		Furthermore, locally around $\bar x$, the mappings 
		$\Phi_\textup{w}+\mathcal C$ and $\Sigma+\mathcal C$ are identical by the weak domination property.
		Noting that the coderivative is a local construction, the claim follows.
\end{proof}

We would like to point the reader's attention to the similarities in the estimates
\eqref{eq:upper_estimate_frontier_map} and \eqref{eq:upper_estimate_weak_frontier_map}.
Indeed, the right-hand side is the same in both estimates which is somewhat unsatisfying
as the precise underlying notion of nondominance seems to make no difference.
	More precisely, it does not seem to be too reasonable to try to find assumptions which
	guarantee correctness of the results, e.g., in \cite{HuyMordukhovichYao2008} for the
	frontier mapping as we obtained the same upper estimate for the coderivative of the
	weak frontier mapping under most likely milder conditions.

With the aid of \cref{lem:coderivative_estimates}, we are in position to obtain an upper estimate
for the coderivative of $\Phi_\textup{w}$ in terms of initial problem data.
The following result is inspired by \cite[Theorem~3.3]{HuyMordukhovichYao2008}.

\begin{theorem}\label{thm:coderivative_estimate}
	Let $(\bar x,\bar z)\in\gph\Phi_\textup{w}$ be fixed and let $\gph\Phi_\textup{w}$ be 
	locally closed around $(\bar x,\bar z)$.
	Assume that \eqref{eq:parametric_problem} possesses the weak domination property for all $x\in U$
	where $U\subset\R^n$ is a neighborhood of $\bar x$.
	Furthermore, let $\Gamma$ be locally bounded at $\bar x$.
	Let $f$ be locally Lipschitz continuous at all points $(\bar x,\bar y)$ 
	such that $\bar y\in\Psi_\textup{w}(\bar x)$ and $f(\bar x,\bar y)=\bar z$ hold.
	Then, for each $z^*\in \mathcal{C}^*_>$,
	\[
		D^*\Phi_\textup{w}(\bar x,\bar z)(z^*)
		\subset
		\left\{
			x_1^*+x_2^*\,\middle|\,
			\begin{aligned}
				&\bar y\in\Psi_\textup{w}(\bar x),\,f(\bar x,\bar y)=\bar z,
				\\
				&(x_1^*,y^*)\in D^*f(\bar x,\bar y)(z^*),\,
				x_2^*\in D^*\Gamma(\bar x,\bar y)(y^*)
			\end{aligned}
		\right\}.
	\]
	If, additionally, $f$ is continuously differentiable at all points $(\bar x,\bar y)$ 
	such that $\bar y\in\Psi_\textup{w}(\bar x)$ and $f(\bar x,\bar y)=\bar z$ hold, 
	then, for each $z^*\in \mathcal{C}^*_>$,
	\[
		D^*\Phi_\textup{w}(\bar x,\bar z)(z^*)
		\subset
		\left\{
			f'_x(\bar x,\bar y)^\top z^*+x^*\,\middle|\,
			\begin{aligned}
				&\bar y\in\Psi_\textup{w}(\bar x),\,f(\bar x,\bar y)=\bar z,
				\\
				&x^*\in D^*\Gamma(\bar x,\bar y)\bigl(f'_y(\bar x,\bar y)^\top z^*\bigr)
			\end{aligned}
		\right\}.
	\]
\end{theorem}
\begin{proof}
	Let us note that closedness of $\gph\Gamma$ and local boundedness of $\Gamma$ at $\bar x$ yield that, 
	locally around $(\bar x,\bar z)$, 
	$\gph\Sigma$ and $\gph(\Sigma+\mathcal C)$ are closed.
	Combining \cref{lem:coderivative_estimates} and the first part of the proof of \cite[Theorem~3.3]{HuyMordukhovichYao2008}
	while taking into account that $-N_{\mathcal C}(0)=\mathcal{C}^*\supset \mathcal{C}^*_>$ since $\mathcal C$ is a cone gives
	\[
		D^*\Phi_\textup{w}(\bar x,\bar z)(z^*)\subset D^*\Sigma(\bar x,\bar z)(z^*)
	\]
	for each $z^*\in \mathcal{C}^*_>$.
	
	Defining $\widetilde\Gamma\colon\R^n\tto\R^n\times\R^m$ by $\widetilde\Gamma(x):=\{x\}\times\Gamma(x)$ for each 
	$x\in\R^n$ gives $\Sigma=f\circ\widetilde\Gamma$, and
	\[
		D^*\widetilde\Gamma(\bar x,(\bar x,y))(x^*,y^*)
		=
		\{x^*\}+D^*\Gamma(\bar x,y)(y^*)
	\]
	holds for each $y\in\Gamma(\bar x)$ and $(x^*,y^*)\in\R^n\times\R^m$, see e.g.\
	\cite[Lemma~5.7]{BenkoMehlitz2022}.
	Thus, in order to infer the more general assertion of the theorem, 
	we apply the chain rule from \cite[Theorem~5.2]{BenkoMehlitz2022}.
	Therefore, we note that the involved intermediate mapping
	\begin{equation}\label{eq:intermediate_map}
		(x,z)\mapsto \{(x,y)\,|\,y\in\Gamma(x),\,f(x,y)=z\}
	\end{equation}
	is inner semicompact at $(\bar x,\bar z)$ w.r.t.\ its domain
	since $\Gamma$ is assumed to be locally bounded at $\bar x$, 
	and that $f$ is Lipschitz continuous at $(\bar x,\bar y)$ which serves as a 
	qualification condition. 
	Then we find
	\begin{align*}
		D^*\Sigma(\bar x,\bar z)(z^*)
		\subset
		\left\{
			\tilde x^*
			\,\middle|\,
			\begin{aligned}
				&\bar y\in\Gamma(\bar x),\,f(\bar x,\bar y)=\bar z,
				\\
				&(x^*,y^*)\in D^*f(\bar x,\bar y)(z^*),\,
				\tilde x^*\in D^*\widetilde\Gamma(\bar x,(\bar x,\bar y))(x^*,y^*)
			\end{aligned}
		\right\},
	\end{align*}
	and noting that,
	due to \cref{prop:value_function_char_of_eff_set}, 
	each $\bar y\in\Gamma(\bar x)$ with $f(\bar x,\bar y)=\bar z$ already belongs to
	$\Psi_\textup{w}(\bar x)$ yields the claim.
	
	The second statement of the theorem is an obvious consequence of the first assertion.
\end{proof}

Let us note that the chain rule from \cite[Theorem~5.2]{BenkoMehlitz2022} also allows 
for the derivation of refined upper estimates in terms of directional coderivatives.

We would like to point the interested reader to the fact that in the proofs of
\cite[Theorem~3.3]{HuyMordukhovichYao2008} and \cite[Theorem~4.1]{LafhimZemkoho2022}, the authors
rely on the chain rule from \cite[Theorem~3.18(i)]{Mordukhovich2006} which is only possible if
the intermediate mapping from \eqref{eq:intermediate_map} possesses certain inner semicontinuity
properties. These, however, are not inherent from the assumptions stated in \cite{HuyMordukhovichYao2008,LafhimZemkoho2022},
and one can easily check that even inner semicontinuity of $\Gamma$ is not sufficient for
these requirements.
Thus, the aforementioned results are not reliable.
Let us also note that, in situations where $f$ is Lipschitzian, 
the coderivative of $\Sigma$ can be computed exploiting a combination of the
image and pre-image rule as we have
\[
	\gph\Sigma = \Pi(\{(x,y,z)\,|\,(x,y,f(x,y)-z)\in\gph\Gamma\times\{0\}\})
\]
where $\Pi\colon\R^n\times\R^m\times\R^q\to\R^n\times\R^q$ is given by $\Pi(x,y,z):=(x,z)$ for all
$x\in\R^n$, $y\in\R^m$, and $z\in\R^q$, see e.g.\ \cite[Sections~5.1.3, 5.1.4]{BenkoMehlitz2022}.
This, however, still makes some additional boundedness assumption on $\Gamma$ necessary, and the resulting
upper estimate on the coderivative of $\Sigma$ would also comprise a union over $\Psi_\textup{w}(\bar x)$.
More precisely, this approach precisely recovers the estimates from
\cref{thm:coderivative_estimate}.

Let us specify \cref{thm:coderivative_estimate} for the scalar case where $q:=1$ and $\mathcal C:=\R_+$.
Then the weak domination property is inherently satisfied.
Furthermore, we find $1\in \mathcal{C}^*_>=(0,\infty)$, 
and $D^*\Phi_\textup{w}(\bar x,\varphi(\bar x))(1)$ equals the so-called limiting subdifferential
at $\bar x$ 
of the optimal valued function $\varphi$ from \eqref{eq:optimal_value_function}
provided the latter is continuous around $\bar x$, see \cite[Definition~1.77, Theorem~1.80]{Mordukhovich2006}.
Then the estimates in \cref{thm:coderivative_estimate} recover those obtained
in \cite[Theorem~7]{MordukhovichNamYen2009} under slightly different assumptions.
Hence, one can interpret \cref{thm:coderivative_estimate} as a generalization of this result
to the multiobjective situation.
Nevertheless, \cref{thm:coderivative_estimate} is of limited practical use as it provides no information about
situations where the argument $z^*$ of the coderivative is chosen in $\R^q\setminus \mathcal{C}^*_>$, and such
values of $z^*$ might be of interest when deriving necessary optimality conditions for
\eqref{eq:VFR_weff} as the concept of subdifferentials is not reasonable in vector-valued scenarios
whenever scalarization is not taken into account.

\subsection{Coderivative estimates via scalarization}\label{sec:coderivative_scalarization}

In this section, we aim to find upper coderivative estimates for the weak frontier mapping $\Phi_\textup{w}$ via
a scalarization approach. This, however, makes additional assumptions necessary.
\begin{assumption}\label{ass:scalarization_approach}
	For each $x\in\R^n$, $f(x,\cdot)$ is $\mathcal C$-convex and $\Gamma(x)$ is convex.
\end{assumption}

We emphasize that \cref{ass:scalarization_approach} guarantees that $\Sigma(x)+\mathcal C$ is convex for each $x\in\R^n$.
Hence, due to \cite[Corollary~5.29]{Jahn2011}, for each $\bar x\in\dom\Psi_\textup{w}$, we find the equivalence
\[
	y\in\Psi_\textup{w}(\bar x)
	\quad\Longleftrightarrow\quad
	\exists \lambda\in \mathcal{C}^*\cap\mathbb S_1(0)\colon\quad
	y\in\Psi_\textup{w}^\textup{sc}(\bar x,\lambda)
\] 
where $\Psi_\textup{w}^\textup{sc}\colon\R^n\times\R^q\tto\R^m$ is the solution mapping associated with the scalar
parametric optimization problem
\begin{equation}\label{eq:scalarized_parametric_optimization_problem}\tag{P$(x,\lambda)$}
	\min\limits_y\{\langle \lambda,f(x,y)\rangle\,|\, y\in\Gamma(x)\}.
\end{equation}

Based on the chain rule, we obtain a simple upper estimate for the coderivative of $\Phi_\textup{w}$ in terms
of the coderivative of $\Psi_\textup{w}$, as we have $\Phi_\textup{w}=f\circ\widetilde\Psi_\textup{w}$ where
$\widetilde\Psi_\textup{w}\colon\R^n\tto\R^n\times\R^m$ is given by 
$\widetilde\Psi_\textup{w}(x):=\{x\}\times\Psi_\textup{w}(x)$ for each $x\in\R^n$. 
We omit the proof as it mainly parallels the arguments used to
validate \cref{thm:coderivative_estimate}.
\begin{lemma}\label{lem:upper_estimate_via_Coderivative_of_solutions}
	Fix $(\bar x,\bar z)\in\gph\Phi_\textup{w}$ and let $\gph\Phi_\textup{w}$ be locally closed around $(\bar x,\bar z)$.
	Let $\Psi_\textup{w}$ be inner semicompact at $\bar x$ w.r.t.\ $\dom\Psi_\textup{w}$.
	Furthermore, for each $\bar y\in\Psi_\textup{w}(\bar x)$ with $f(\bar x,\bar y)=\bar z$, 
	let $\gph\Psi_\textup{w}$ be closed locally around $(\bar x,\bar y)$,
	and let $f$ be locally Lipschitz continuous at $(\bar x,\bar y)$.
	Then, for each $z^*\in\R^q$,
	\[
		D^*\Phi_\textup{w}(\bar x,\bar z)(z^*)
		\subset
		\left\{
			x_1^*+x_2^*\,\middle|\,
			\begin{aligned}
				&\bar y\in\Psi_\textup{w}(\bar x),\,f(\bar x,\bar y)=\bar z,
				\\
				&(x_1^*,y^*)\in D^*f(\bar x,\bar y)(z^*),\,
				x_2^*\in D^*\Psi_\textup{w}(\bar x,\bar y)(y^*)
			\end{aligned}
		\right\}.
	\]
	If, additionally, $f$ is continuously differentiable at all points $(\bar x,\bar y)$ 
	such that $\bar y\in\Psi_\textup{w}(\bar x)$ and $f(\bar x,\bar y)=\bar z$ hold, 
	then, for each $z^*\in \R^q$,
	\[
		D^*\Phi_\textup{w}(\bar x,\bar z)(z^*)
		\subset
		\left\{
			f'_x(\bar x,\bar y)^\top z^*+x^*\,\middle|\,
			\begin{aligned}
				&\bar y\in\Psi_\textup{w}(\bar x),\,f(\bar x,\bar y)=\bar z,
				\\
				&x^*\in D^*\Psi_\textup{w}(\bar x,\bar y)\bigl(f'_y(\bar x,\bar y)^\top z^*\bigr)
			\end{aligned}
		\right\}.
	\]
\end{lemma}

Let us emphasize that \cref{lem:upper_estimate_via_Coderivative_of_solutions} is correct even in the
absence of \cref{ass:scalarization_approach} which is not needed for the proof.
We note that the estimates from \cref{lem:upper_estimate_via_Coderivative_of_solutions} are closely related to
those ones obtained in \cref{thm:coderivative_estimate}.
Indeed, the only difference (apart from the underlying assumptions) is that in \cref{thm:coderivative_estimate},
all estimates are valid for $z^*\in \mathcal{C}^*_>$ only, while the coderivative of the efficiency mapping $\Psi_\textup{w}$
appearing in \cref{lem:upper_estimate_via_Coderivative_of_solutions} is replaced by the coderivative of the
feasibility mapping $\Gamma$ in \cref{thm:coderivative_estimate}.
Let us also visualize this in exemplary way.

One can easily check that the estimate provided in \cref{lem:upper_estimate_via_Coderivative_of_solutions}
is actually sharp in the context of \cref{ex:coderivative_of_weak_frontier_map}.
Indeed, one obtains
\[
	D^*\Phi_\textup{w}(\bar x,\bar z)(z^*)
	=
	D^*\Psi_\textup{w}(\bar x,\bar y)(z_2^*)
	=
	\begin{cases}
		\{0\}	&	z_2^*\neq 0,\\
		\R_+	&	z_2^*=0
	\end{cases}
\]
for each $z^*\in\R^2$ in this situation 
(observe that the only point $y\in\Psi_\textup{w}(\bar x)$ which satisfies  $f(\bar x,y)=\bar z$ 
in the present situation is $\bar y$). 
We also note that the estimates from \cref{thm:coderivative_estimate} do not apply for
$z^*\in\R\times\{0\}$ as these vectors do not belong to $\intr\R^2_+$.

It is an essential disadvantage of the estimates from \cref{lem:upper_estimate_via_Coderivative_of_solutions} 
that they are not stated in terms of initial problem data as the coderivative of $\Psi_\textup{w}$ is a
fully implicit object. Subsequently, we aim to estimate it from above. In order to deal with this issue,
\cref{ass:scalarization_approach} turns out to be helpful as we will demonstrate below.
For brevity of notation, we introduce an intermediate mapping $\Xi\colon\R^n\times\R^m\tto\R^q$ given by
\[
	\forall x\in\R^n\,\forall y\in\R^m\colon\quad
	\Xi(x,y):=\{\lambda\in \mathcal{C}^*\cap\mathbb S_1(0)\,|\,y\in\Psi_\textup{w}^\textup{sc}(x,\lambda)\}.
\]
For given $(x,y)\in\gph\Psi_\textup{w}$, $\Xi(x,y)$ comprises all (normalized) scalarization parameters 
for which $y$ as a minimizer of \eqref{eq:scalarized_parametric_optimization_problem}.
We note that, by compactness of $\mathcal{C}^*\cap \mathbb S_1(0)$, $\Xi$ is inner semicompact w.r.t.\ its domain
at each point of its domain.

\begin{lemma}\label{lem:coderivative_efficiency_via_Coderivative_scalarized}
	Fix $(\bar x,\bar y)\in\gph\Psi_\textup{w}$ and let $\gph\Psi_\textup{w}$ be closed locally around $(\bar x,\bar y)$.
	Furthermore, for each $\lambda\in\Xi(\bar x,\bar y)$, let $\gph\Psi_\textup{w}^\textup{sc}$ be closed
	locally around $((\bar x,\lambda),\bar y)$.
		Additionally, let \cref{ass:scalarization_approach} hold.
		Finally, let one of the following conditions be valid:
	\begin{enumerate}
		\item\label{item:polyhedrality} 
			$\Psi_\textup{w}^\textup{sc}$ is a polyhedral set-valued mapping,
			$\mathcal C$ is a polyhedral cone, and the unit sphere $\mathbb S_1(0)$ is taken w.r.t.\ the $1$- or $\infty$-norm,
		\item\label{item:Aubin_property} 
			for each $\lambda\in\Xi(\bar x,\bar y)$, the qualification condition
			\[
				(0,\eta)\in D^*\Psi_\textup{w}^\textup{sc}((\bar x,\lambda),\bar y)(0),\,
				\eta\in -N_{\mathcal{C}^*\cap\mathbb S_1(0)}(\lambda)
				\quad\Longrightarrow\quad
				\eta=0
			\]
			holds.
	\end{enumerate}
	Then, for each $y^*\in\R^m$, we have
	\[
		D^*\Psi_\textup{w}(\bar x,\bar y)(y^*)
		\subset
		\left\{
			x^*\,\middle|\,
			\begin{aligned}
				&(x^*,\eta)\in D^*\Psi_\textup{w}^\textup{sc}((\bar x,\lambda),\bar y)(y^*),\\
				&\lambda\in\Xi(\bar x,\bar y),\,
				\eta\in -N_{\mathcal{C}^*\cap \mathbb S_1(0)}(\lambda)
			\end{aligned}
		\right\}.
	\]
\end{lemma}
\begin{proof}
	Defining $\Upsilon\colon\R^n\tto\R^n\times\R^q$ by $\Upsilon(x):=\{x\}\times(\mathcal{C}^*\cap\mathbb S_1(0))$ for each $x\in\R^n$,
	we find $\Psi_\textup{w}=\Psi_\textup{w}^\textup{sc}\circ\Upsilon$ since \cref{ass:scalarization_approach} holds.
	Since $\Xi$ is inner semicompact at $(\bar x,\bar y)$ w.r.t.\ its domain, the intermediate mapping
	\[
		(x,y)
		\mapsto
		\{(x,\lambda)\in\Upsilon(x)\,|\,y\in\Psi^\textup{sc}_\textup{w}(x,\lambda)\}
		=
		\{(x,\lambda)\,|\,\lambda\in\Xi(x,y)\}
	\]
	is inner semicompact at $(\bar x,\bar y)$ w.r.t.\ its domain as well.
	Thus, the result follows from the chain rule stated in \cite[Theorem~5.2]{BenkoMehlitz2022} while observing that
	\[
		D^*\Upsilon(\bar x,(\bar x,\lambda))(x^*,\eta)
		=
		\begin{cases}	
			\{x^*\}	&	\eta\in -N_{\mathcal{C}^*\cap\mathbb S_1(0)}(\lambda),\\
			\emptyset	&	\text{otherwise}
		\end{cases}
	\]
	holds for each $\lambda\in\Xi(\bar x,\bar y)$ and $(x^*,\eta)\in\R^n\times\R^q$ while each of the additionally postulated
	assumptions guarantees validity of the qualification condition which is imposed in \cite[Theorem~5.2]{BenkoMehlitz2022}.
\end{proof}

Let us note that \cref{lem:coderivative_efficiency_via_Coderivative_scalarized} generalizes 
\cite[Proposition~3.2]{Zemkoho2016} to arbitrary ordering cones while coming along with a simpler proof and another
type of qualification condition.
Furthermore, we would like to mention 
that \cref{lem:coderivative_efficiency_via_Coderivative_scalarized}
can be used to derive necessary optimality conditions for
\eqref{eq:BPP_weff} directly.

In order to obtain a valuable coderivative estimate for the weak frontier mapping from 
\cref{lem:upper_estimate_via_Coderivative_of_solutions,lem:coderivative_efficiency_via_Coderivative_scalarized},
it remains to estimate the appearing coderivative of $\Psi_\textup{w}^\textup{sc}$.
This, however, is a well-established subject in variational analysis since 
$\Psi_\textup{w}^\textup{sc}$ is the solution mapping of a scalar parametric optimization problem.
Exemplary, let us mention that, in the presence of \cref{ass:scalarization_approach}, 
\eqref{eq:scalarized_parametric_optimization_problem} is a convex optimization problem for each
$x\in\R^n$ and $\lambda\in \mathcal{C}^*$. Assuming, additionally, that $f$ is a smooth function, we find
\[
	\forall x\in\R^n\,\forall \lambda\in \mathcal{C}^*\colon\quad
	\Psi_\textup{w}^\textup{sc}(x,\lambda)
	=
	\left\{
		y
		\,\middle|\,
		-f'_y(x,y)^\top \lambda\in N_{\Gamma(x)}(y)
	\right\},
\]
i.e., $\Psi_\textup{w}^\textup{sc}$ possesses the structure 
of a so-called parametric variational system.
Consulting e.g.\ \cite[Section~4.4]{Mordukhovich2006}, coderivative estimates for $\Psi_\textup{w}^\textup{sc}$
can be obtained in terms of derivatives of $f$ and coderivatives of the normal cone mapping
$(x,y)\mapsto N_{\Gamma(x)}(y)$ under suitable constraint qualifications.
The remaining coderivative of the normal cone map may be specified for diverse different types of constraint
systems, see e.g.\ \cite{BenkoGfrererOutrata2000,GfrererOutrata2016,GfrererOutrata2016b,GfrererOutrata2017}.
Exemplary, this has been worked out in the proof of \cite[Corollary~3.4]{Zemkoho2016} in the particular situation
where $\mathcal C:=\R^q_+$ and $\Gamma$ is induced by smooth parametric inequality constraints.

As the purpose of this note is mainly focused on the value behind the model reformulations 
\eqref{eq:VFR_eff} and \eqref{eq:VFR_weff}, we abstain from the postulation of any more technical results, but
simply combine 
\cref{lem:upper_estimate_via_Coderivative_of_solutions,lem:coderivative_efficiency_via_Coderivative_scalarized}
to obtain a coderivative estimate which can be made fully explicit by respecting the above comments by the 
interested reader.

\begin{theorem}\label{thm:coderivative_estimate_scalarization}
	Fix $(\bar x,\bar z)\in\gph\Phi_\textup{w}$ and let $\gph\Phi_\textup{w}$ be locally closed around $(\bar x,\bar z)$.
	Furthermore, let \cref{ass:scalarization_approach} be valid.
		Let $\Psi_\textup{w}$ be inner semicompact at $\bar x$ w.r.t.\ $\dom\Psi_\textup{w}$.
			Furthermore, for each $\bar y\in\Psi_\textup{w}(\bar x)$ with $f(\bar x,\bar y)=\bar z$, 
			let $\gph\Psi_\textup{w}$ be closed locally around $(\bar x,\bar y)$,
			let $f$ be locally Lipschitz continuous at $(\bar x,\bar y)$, 
			and, for each $\lambda\in\Xi(\bar x,\bar y)$, let $\gph\Psi_\textup{w}^\textup{sc}$ be closed
			locally around $((\bar x,\lambda),\bar y)$, and
			let one of the conditions~\ref{item:polyhedrality} or~\ref{item:Aubin_property} of
			\cref{lem:coderivative_efficiency_via_Coderivative_scalarized} be valid.
			Then, for each $z^*\in\R^q$,
				\[
				D^*\Phi_\textup{w}(\bar x,\bar z)(z^*)
				\subset
				\left\{
					x_1^*+x_2^*\,\middle|\,
					\begin{aligned}
						&\bar y\in\Psi_\textup{w}(\bar x),\,f(\bar x,\bar y)=\bar z,\\
						&(x_1^*,y^*)\in D^*f(\bar x,\bar y)(z^*),\\
						&(x_2^*,\eta)\in D^*\Psi_\textup{w}^\textup{sc}((\bar x,\lambda),\bar y)(y^*),\\
						&\lambda\in\Xi(\bar x,\bar y),\,\eta\in -N_{\mathcal{C}^*\cap\mathbb S_1(0)}(\lambda)
					\end{aligned}
				\right\}.
			\]
			If, additionally, $f$ is continuously differentiable at all points $(\bar x,\bar y)$ 
			such that $\bar y\in\Psi_\textup{w}(\bar x)$ and $f(\bar x,\bar y)=\bar z$ hold, 
			then, for each $z^*\in \R^q$,
			\[
				D^*\Phi_\textup{w}(\bar x,\bar z)(z^*)
				\subset
				\left\{
					f'_x(\bar x,\bar y)^\top z^*+x^*\,\middle|\,
					\begin{aligned}
						&\bar y\in\Psi_\textup{w}(\bar x),\,f(\bar x,\bar y)=\bar z,\\
						&(x^*,\eta)\in D^*\Psi_\textup{w}^\textup{sc}((\bar x,\lambda),\bar y)\bigl(f'_y(\bar x,\bar y)^\top z^*\bigr),\\
						&\lambda\in\Xi(\bar x,\bar y),\,\eta\in -N_{\mathcal{C}^*\cap\mathbb S_1(0)}(\lambda)
					\end{aligned}
				\right\}.
			\]
\end{theorem}

Let us note that, based on \cref{thm:closedness_of_graphs}, one can easily check that the numerous closedness
assumptions as well as the inner semicompactness assumption on $\Psi_\textup{w}$ in \cref{thm:coderivative_estimate_scalarization}
are inherently satisfied whenever $\Gamma$ is lower semicontinuous w.r.t.\ its domain locally around $\bar x$, and
locally bounded at $\bar x$. Both assumptions are inherent whenever $\Gamma(x):=\Omega$ holds for some nonempty and compact
set $\Omega\subset\R^m$ and all $x\in\R^n$.
Then it remains to validate condition~\ref{item:polyhedrality} or~\ref{item:Aubin_property} 
of \cref{lem:coderivative_efficiency_via_Coderivative_scalarized}. 

\begin{remark}\label{rem:scalarization_approach_via_value_function}
	Let \cref{ass:scalarization_approach} be valid
	and $\varphi^\textup{sc}\colon\R^n\times\R^q\to\overline\R$ be the optimal value function 
	associated with \eqref{eq:scalarized_parametric_optimization_problem}, i.e.,
	\[
		\forall x\in\R^n\,\forall\lambda\in\R^q\colon\quad
		\varphi^\textup{sc}(x,\lambda)
		:=
		\inf\limits_y\{\langle\lambda,f(x,y)\rangle\,|\,y\in\Gamma(x)\}.
	\]
	Then we find the representation
	\[
		\forall x\in\R^n\colon\quad
		\Phi_\textup{w}(x)
		=
		\{
			z\in\Sigma(x)
			\,|\,
			\exists \lambda\in \mathcal{C}^*\cap\mathbb S_1(0)\colon\,
				\langle\lambda,z\rangle=\varphi^\textup{sc}(x,\lambda)
		\},
	\]
	where $\Sigma\colon\R^n\tto\R^q$ has been defined in \eqref{eq:surrogate_maps_coderivative_calculus}.
	Using $\Upsilon\colon\R^n\tto\R^n\times\R^q$ from the proof of 
	\cref{lem:coderivative_efficiency_via_Coderivative_scalarized} as well as
	$\Upsilon_1\colon\R^n\times\R^q\to\R^q\times\R$ and $\Upsilon_2\colon\R^q\times\R\to\R^q$ given by
	\[
		\forall x\in\R^n\,\forall\lambda\in\R^q\,\forall\alpha\in\R\colon\quad
		\Upsilon_1(x,\lambda):=(\lambda,\varphi^\textup{sc}(x,\lambda)),\qquad
		\Upsilon_2(\lambda,\alpha):=\{z\,|\,\langle\lambda,z\rangle-\alpha=0\},
	\]
	we find $\Phi_\textup{w}(x)=\Sigma(x)\cap(\Upsilon_2\circ(\Upsilon_1\circ\Upsilon))(x)$.
	Seemingly, this representation offers the possibility to obtain an upper estimate for the
	coderivative of $\Phi_\textup{w}$ via an intersection rule, 
	see e.g.\ \cite[Proposition~3.20]{Mordukhovich2006}, and the chain rule,
	see e.g.\ \cite[Theorem~5.2]{BenkoMehlitz2022}.
	Indeed, in the case where $\varphi^\textup{sc}$ is a locally Lipschitzian function, 
	which holds under fairly mild assumptions, one can apply the
	chain rule from \cite{BenkoMehlitz2022} twice in order to find an upper estimate of
	the coderivative of $\Upsilon_2\circ(\Upsilon_1\circ\Upsilon)$ which comprises
	the limiting subdifferential of (a non-specified multiple of) $\varphi^\textup{sc}$.
	The latter, however, cannot be approximated from above by the Cartesian product
	of the associated subdifferentials w.r.t.\ the variable $x$ and $\lambda$ in general
	which is a first difficulty.
	A second problem pops up as the constraint qualification from \cite[Proposition~3.20]{Mordukhovich2006}
	which is sufficient to get the intersection rule working is also likely to fail.
	Hence, we abstain here from providing the details as this approach does not seem too reasonable.
\end{remark}

	\begin{remark}\label{rem:reasonability}
		In scalar bilevel optimization, 
		the value function reformulation is often preferred over the formulation
		comprising the lower-level solution mapping when the derivation of necessary
		optimality conditions is considered.
		This is due to the fact that the computation of upper
		estimates for generalized derivatives of the lower-level value function is
		often much easier than obtaining upper estimates for generalized derivatives 
		of the lower-level solution mapping. 
		Thus, on the one hand,  
		\cref{lem:upper_estimate_via_Coderivative_of_solutions} 
		and 
		\cref{thm:coderivative_estimate_scalarization} 
		seemingly are of limited practical use 
		in the context of necessary optimality conditions as the coderivative 
		of some solution mapping associated with the lower-level problem appears.
		Indeed, one could try to directly make use of 
		\cref{lem:coderivative_efficiency_via_Coderivative_scalarized}
		in order to obtain necessary optimality conditions.
		On the other hand, the value function approach to bilevel optimization is
		often preferred for the construction of solution algorithms, 
		and taking this aspect into account, some upper estimate for the
		coderivative of the weak frontier mapping might be essential, and this is
		where \cref{lem:upper_estimate_via_Coderivative_of_solutions} 
		and 
		\cref{thm:coderivative_estimate_scalarization} 
		may become handy.
	\end{remark}

\section{Conclusions}\label{sec:conclusions}

In this paper, we visualized that multiobjective bilevel optimization is, up to a certain degree,
only reasonable if the lower-level efficiency-type and frontier-type mapping possess closed
graphs as this is a fundamental assumption in order to guarantee the existence of solutions (in
whatever sense) and to allow for the derivation of optimality conditions.
We presented conditions which guarantee validity of these closedness properties, gave existence results, and
discussed the computation of upper coderivative estimates for frontier-type mappings.
The latter is an essential step on the route to necessary optimality conditions as visualized in the recent
paper \cite{LafhimZemkoho2022} which, however, has some flaws which we carved out here.
Furthermore, we commented on different ways on how to introduce a value function reformulation of
a multiobjective bilevel optimization problem which might be important in algorithmic applications.

In the future, it has to be investigated whether the value function reformulation of a multiobjective
optimization problem can be exploited beneficially for its algorithmic solution.
Whenever the lower-level problem is a scalar optimization problem, this has already been demonstrated
in several different ways in the literature, see e.g.\ \cite[Section~20.6.4]{Dempe2020} for an overview.
	Furthermore, it remains open whether the approach from \cref{sec:intermediate_approach}
	is suitable from the viewpoint of optimality conditions as, in contrast to 
	\cref{sec:coderivative_calculus}, it is not clear how the derivation of coderivative
	estimates for the mappings $\overline\Phi$ and $\overline\Psi$ can be done.


\end{document}